\documentclass[11pt,twoside]{amsart}
\usepackage{amsfonts}
\usepackage{epsfig,graphics}
\usepackage{amssymb}
\usepackage{amscd}
\usepackage{accents}
\usepackage{stmaryrd}

\newtheorem{theorem}{Theorem}[section]

\newtheorem{proposition}[theorem]{Proposition}
       
\newtheorem{definition}[theorem]{Definition}

\newtheorem{corollary}[theorem]{Corollary}

\newtheorem{remark}[theorem]{Remark}

\newcommand{\Sym}{\text{Sym}}

\newcommand{\Ad}{\text{Ad }}
\newcommand{\End}{\text{End}}

\newcommand{\im}{\text{im }}

\newcommand{\GL}{\text{GL}}

\newcommand{\id}{\text{Id}}

\newcommand{\BZ}{{\bf Z}}
\newcommand{\BN}{{\bf N}} 
\newcommand{\BR}{{\bf R}}

\newcommand{\alggl}{{\mathfrak{g}\mathfrak{l}}}

\newcommand{\intnu}{\accentset{\circ}{\nu}}

\newenvironment{Pf}{\medskip \noindent {\bf Proof: }}
   {$\diamondsuit$ }

\begin{document}
\pagenumbering{arabic}
\title[Nonstationary geometric structures]{Nonstationary smooth geometric structures for contracting measurable cocycles}
\author{Karin Melnick}
\thanks{This article was in large part written during a visit to the Max Planck Institut f\"ur Mathematik in Bonn, whose support I gratefully acknowledge, along with partial support from NSF grant DMS 1255462.  I wish to acknowledge the important role of David Fisher in this project, as the one who proposed using the results of Ruelle for the measurable setting, and who also introduced me to Feres' paper.  I furthermore thank him for numerous helpful discussions and insights.} 

\date{\today}
\maketitle

\begin{abstract}
We implement a differential-geometric approach to normal forms for contracting measurable cocycles to $\mbox{Diff}^q(\BR^n, {\bf 0})$, $q \geq 2$.  We obtain resonance polynomial normal forms for the contracting cocycle and its centralizer, via $C^q$ changes of coordinates.  These are interpreted as nonstationary invariant differential-geometric structures.  We also consider the case of contracted foliations in a manifold, and obtain $C^q$ homogeneous structures on leaves for an action of the group of subresonance polynomial diffeomorphisms together with translations. 
\end{abstract}

\section{Introduction}

Normal forms for differentiable contractions have a long history.  Sternberg proved in \cite{sternberg.normal.forms} that a $C^q$-smooth diffeomorphism $f$ of $\BR^n$ fixing ${\bf 0}$ and satisfying $\Vert D_{\bf 0} f \Vert < 1$ is $C^q$ conjugate on a neighborhood of ${\bf 0}$ to a polynomial of \emph{resonance} type, obtained from the Taylor series of $f$ at ${\bf 0}$.  The degree of this polynomial is bounded in terms of the spectrum of $D_{\bf 0} f$.  

This article concerns a generalization of this normal forms theorem to a flow acting on a product $M \times \BR^n$, where $M$ is a probability space, via a measurable cocycle to $\mbox{Diff}^q(\BR^n,{\bf 0})$.

\subsection{Previous results on nonstationary normal forms}

Guysinsky and Katok found nonstationary normal forms in the $C^0$ setting in \cite{guysinsky.katok.normal} and \cite{guysinsky.nonstationary.etds}.  They consider a transformation $F$ of a continuous $\BR^n$-bundle $\mathcal{E} \rightarrow M$, where $M$ is a compact metric space, acting along fibers by $C^q$ diffeomorphisms preserving the zero section. The contracting assumption is that $\sup_{x \in M} \Vert D_{\bf 0} F_x \Vert < 1$.  Assuming that the Mather spectrum of $F$ acting on $\Gamma(\mathcal{E})$ has the \emph{narrow bands} property, they prove that a continuous family of $C^q$ coordinate changes leaves $F$ acting along fibers by \emph{subresonance polynomials}.  Their degree is bounded in terms of the Mather spectrum.  

The outline of their proof resembles Sternberg's.  They first solve formally and then make a fixed point argument on the space of coordinate changes to obtain the solution.  An important element in their conclusions is that the centralizer of $F$ also acts along fibers by subresonance polynomials in the chosen coordinates.   

In \cite{feres.normal.forms}, Feres gives an interpretation of Guysinksy and Katok's normal forms in the case $\mathcal{E}$ represents an $F$-invariant lamination in $M$, as a continuous $F$-invariant family of geometric structures: they are reductions of the frame bundles of some sufficiently high order of the leaves to the group of subresonance polynomials.  He further obtains \emph{generalized connections} on these reductions, from which he constructs the polynomial charts on leaves using Gromov's Frobenius theorem for partial differential relations---all under the additional assumption that the Mather spectral partition is differentiable along leaves.

This differentiability assumption on the spectral partition is problematic, because it does not usually hold.  In a recent paper, Kalinin and Sadovskaya show \cite{kalinin.sadovskaya.normal.smoothness} that differentiability along leaves of the Mather spectral \emph{filtration}, which does hold in general, suffices to obtain smoothness along leaves of the normal forms.

These $C^0$ nonstationary normal forms have been used to prove numerous outstanding theorems, primarily on local rigidity of discrete group actions on homogeneous spaces (to cite just a few: \cite{fisher.margulis.inv}, \cite{fisher.margulis.annals}, \cite{katok.spatzier.locrigidity}, \cite{kalinin.sadovskaya.tns.anosov.Zk}).  The initial action satisfies the narrow band condition on the spectrum, and sufficiently small deformations retain the property.  We are interested in other situations in which one cannot necessarily verify the narrow band condition.  In this case, an alternative approach is to consider measurable contractions.

\subsection{The approach in this paper}

We construct polynomial normal forms for contracting measurable cocycles to $\mbox{Diff}^q (\BR^n, {\bf 0})$, $q \geq 2$.  Contracting in this context means negative Lyapunov spectrum.  Our primary tool is the theorem of Ruelle in \cite{ruelle.diff.ergodic}, rather than an explicit fixed point argument.  We take a differential-geometric point of view and apply Ruelle's theorem to a prolongation of the initial cocycle to higher-order frame bundles $\mathcal{E}^{(r)}$ of $\mathcal{E}$.  We then find a measurable invariant family of smooth reductions of $\mathcal{E}^{(r)}$ to the group of subresonance polynomials.  From these reductions, at sufficiently high order, we construct, by methods similar to Feres', coordinate atlases on fibers in which our cocyle and its centralizer have values in \emph{resonance polynomials}; these are defined in terms of the Lyapunov spectrum. 
  
A statement of resonance polynomial normal forms for measurable contracting cocycles and their centralizers appears in \cite[Thms 6.1, 6.2]{katok.kalinin.normal}.  They briefly sketch the steps of a proof, following the scheme of \cite{guysinsky.katok.normal}.  A detailed proof for $q = \infty$ was given in 2005 by W. Li and K. Lu \cite{li.lu.nonuniform.normal}.  They in fact treat \emph{hyperbolic} measurable cocycles---that is, they assume only that $0$ is not in the spectrum; they do not, however, extend their polynomial normal forms to the centralizer of the initial transformation.  

Our preprint appeared online essentially simultaneously with one by Kalinin and Sadovskaya \cite{kalinin.sadovskaya.nonuniform.normal} in which they also obtain, by a different approach from ours, \emph{subresonance} polynomial normal forms for contracting measurable cocycles and their centralizers for finite $q$.  They can treat $q = r + \alpha$ for integral $r \geq 1$ and $0 < \alpha < 1$; they require $q + \alpha$ greater than the ratio of minimum and maximum Lyapunov exponents, while we require $q$ greater than this ratio.  They require a temperedness assumption on the $C^q$ norms of their cocycle, while we, in order to apply Ruelle's theorem, make a stronger assumption of finite first moment.  They apply the normal forms theorem to contracted foliations and obtain a homogeneous structure along leaves, as in our theorem \ref{thm.foliation.atlas} below.

The group of resonance polynomials is significantly smaller than the corresponding subresonance group. On the other hand, we lose the continuity between fibers present in the uniform case.  In the important special case that the fibers are plaques of an invariant $C^0$ foliation, however, we obtain $C^q$ smoothness along leaves, and homogeneity of the atlas.  A recent paper \cite{bufetov.etal.hoelder.oseledets} proves H\"older continuity of Lyapunov decompositions for diffeomorphisms of compact manifolds on sets of arbitrarily large measure, as in Lusin's theorem; as the regularity of our frame bundle reductions and resulting atlases seems to be controlled by regularity of the Lyapunov filtration, one could hope for a similar regularity result for them.

We provide an additional differential-geometric interpretation of the invariant structures on fibers, which generalizes the well-known flat affine connections on fibers in the case of 1/2-pinched spectrum: let $r$ be the maximum degree of the resonance polynomials associated to the Lyapunov spectrum, and assume $r \leq q-1$.  Then in almost every fiber there is a $C^{q-1}$ flag of submanifolds $\nu_x^1 \subset \cdots \subset \nu_x^l$, $l \leq r$, and a $C^{q-2}$ flat connection $\nabla_x^i$ on the normal bundle of $\nu_x^{i}$ in $\nu_x^{i+1}$, for each $i = 1, \ldots, l-1$.  In the case of a foliation, these flags of submanifolds fit together to a $C^{q-1}$ flag of foliations, along with their connections, in almost-every leaf.

\subsection{Statement of main results}

Let $( \{ \varphi^t \}, M, \mu)$ be an ergodic flow preserving a probability measure.  Let $\mathcal{E} \cong M \times \BR^n \rightarrow M$ be a measurable $\BR^n$-bundle to which $\{ \varphi^t \}$ lifts, preserving $M \times \{ {\bf 0} \}$, and acting $C^q$-smoothly on the fibers for $q \geq 2$. Denote by $F(t,x) = F^t_x$ the corresponding cocycle to $\mbox{Diff}^q(\BR^n,{\bf 0})$, and by $T(t,x)= T^t_x$ the cocycle $D_{\bf 0} F^t_x$.  Assume that 
\begin{equation}
\label{eqn.met}
\sup_{0 \leq t \leq 1} \ln^+ \Vert T^t \Vert,  \sup_{0 \leq t \leq 1} \ln^+ \Vert T^{-t} \Vert \in L^1(M,\mu) \tag{MET}
\end{equation}
and let 
$$- \infty < \lambda^{(1)} < \cdots < \lambda^{(s)} < 0 \leq \lambda^{(s+1)} < \cdots < \lambda^{(m)}$$ 
be the Lyapunov spectrum of $T$, given by the Multiplicative Ergodic Theorem.  The cocycle $\{ F^t_x \}$ will be called a \emph{measurable contraction} if $m=s$---that is, if the Lyapunov spectrum is all negative.  

See definition \ref{def.resonance.polys} below for the group of \emph{resonance polynomials} $\mathcal{H}^{(r),0}$ associated to the Lyapunov spectrum.  The \emph{centralizer} $Z(F^t_x)$ of the $\{ \varphi^t \}$ action on $\mathcal{E}$ is given in definition \ref{def.centralizer}.
For $F \in \mbox{Diff}^q(\BR^n)$, write $\Vert F \Vert_q$ for the $C^q$-norm of $F$ restricted to the unit disk $\overline{B(1)}$. 

{\bf Theorem} \ref{thm.main}:
\emph{Assume $\{F^t_x\}$ is a measurable contraction, satisfying}  
$$\sup_{0 \leq t \leq 1} \ln^+ \Vert F^{\pm t}_x \Vert_q \in L^1(M,\mu).$$ 
\emph{Let $r = \lfloor \lambda^{(1)}/\lambda^{(s)} \rfloor$, and assume $r \leq q-1$.  Then there exist the following differential-geometric structures on $\mathcal{E}_x$, for all $x$ in a $\varphi^t$-invariant, $\mu$-conull subset:}
\begin{enumerate}
%% \item A $C^{q-r+1}$ reduction of $\mathcal{E}^{(r-1)}$ to $\mathcal{X}^{(r-1)}$, invariant by $F^t_x$, up to the algebraic hull of $J^{(r-1)}_{\bf 0} F^t_x$, which is contained in $\mathcal{H}^{(r-1),0}$.  If there are no resonances in $\Sigma^0$, then this algebraic hull is in $\GL(n)$.  This reduction is also invariant by any $(\{G^k_x \},\psi) \in Z(F^t_x)$, up to the algebraic hull of $J_{\bf 0}^{(r-1)} G^k_x$ in $\mathcal{H}^{(r-1),0}$.

\item \emph{A family $\mathcal{A}_x$ of $C^q$ charts, with transitions in $\mathcal{H}^{(r),0}$ (restricted to a neighborhood of ${\bf 0}$). The collection $\cup_x \mathcal{A}_x$ is $Z(F^t_x)$-invariant. }

\item \emph{A $C^{q-1}$ filtration $\mathcal{V}^1_x \subset \cdots \subset \mathcal{V}^l_x = T \mathcal{E}_x$, with $l \leq r$, equipped with $C^{q-2}$ flat connections $\nabla_x^i$ on $\mathcal{V}_x^i/\mathcal{V}_x^{i-1}$, $i = 1, \ldots, l$.  The filtrations and connections are invariant by $Z(F^t_x)$. }
\end{enumerate}

Elements of $Z(F^t_x)$ act by resonance polynomials in the coordinates given by (1).  Part (1) above should be compared with \cite[Thm 1.1(v)]{li.lu.nonuniform.normal}, while noting that they assume $q=\infty$, and that they do not extend their result to the centralizer.  The filtration in (2) restricted to the zero section of $\mathcal{E}$ is a conglomeration of the Lyapunov filtration of $\{ T^t_x \}$. 

In the important case that $M$ is a compact manifold and $\{ \varphi^t \}$ a $C^0$ flow, and where the fibers of $\mathcal{E}$ are plaques of a $C^0$ foliation $L$ by $n$-dimensional submanifolds admitting a $\varphi^t$-invariant $C^q$-structure (see section \ref{sec.filtered.foliations} below for details), we have the following smooth geometric structures on the leaves:

{\bf Theorem} (compare Thm \ref{thm.foliation.atlas} below)
\emph{Assume $q \geq \lfloor \lambda^{(1)}/\lambda^{(s)} \rfloor + 1$ and that for fixed $t$, the $q$-jet $J_x^{(q)}( \left. \varphi^t \right|_{L})$ is continuous in $x$.  For all $x$ in a $\varphi^t$-invariant, $\mu$-conull subset, for all $y \in L_x$, there is a family of $C^q$ charts on $L_x$ at $y$, varying measurably between leaves, with the additional properties:} 
\begin{enumerate}
\item \emph{They are global diffeomorphisms $\BR^n \rightarrow L_x$;}
\item  \emph{$\varphi^t$ and its centralizer act by resonance polynomials in these charts;} 
\item  \emph{they make $L_x$ into a homogeneous space for the group generated by translations and subresonance polynomials (see definition \ref{def.resonance.polys}).}
\end{enumerate}
% \begin{enumerate} 
% \item Global diffeomorphisms: For all $x$ in a $\varphi^t$-invariant, $\mu$-conull set, each $\beta \in \mathcal{B}^0(y)$ is a diffeomorphism $\beta : (\BR^n, {\bf 0}) \rightarrow (L_x,y)$.
% \item Homogeneity: Given $z \in L_x$, $\beta \in \mathcal{B}^0(y)$, and $\gamma \in \mathcal{B}^0(z)$, we have $(\gamma \circ \beta^{-1}) \cdot \mathcal{B}^0(y) = \mathcal{B}^0(z)$.
% \item Invariance: Let $\mathcal{B}^0_x = \cup_{y \in L_x} B^0(y)$.  The collection $\cup_x \mathcal{B}_x^0$ is $Z(\varphi^t)$-invariant, and in these coordinates, $Z(\varphi^t)$ acts by resonance polynomials---that is, elements of $\mathcal{H}^{(r),0}$.
% \item Structure group: For $\beta, \gamma \in \mathcal{B}_x^0$, the transition $\gamma^{-1} \circ \beta$ acts on $\BR^n$ by a translation composed with the an element of $\mathcal{H}^{(r)}$.
% \end{enumerate}
\emph{If the spectrum $\Sigma^0$ is $1/2$-pinched, then $L$ carries an invariant family of $C^{q-1}$ flat affine structures, measurable in $x$.}

The Lyapunov-Ruelle foliations tangent to the Lyapunov filtration along $L$ carry the following differential-geometric structure:

{\bf Theorem} \ref{thm.foliation.connxns}:
\emph{The foliation $\cup_x L_x$ contains a $Z(\varphi^t)$-invariant filtered family of $C^{q-1}$ subfoliations} 
$$ L^{i_1} \subset \cdots L^{i_{l}} = L \qquad l \leq r,$$
\emph{in which, for all $1 \leq j \leq l$, the subfoliations $L^{i_{j-1}} \subset L^{i_j}$ carry a $Z(\varphi^t)$-invariant, flat transverse connection, which is measurable in $x$ and $C^{q-2}$ inside almost every $L_x$; in particular, almost every leaf $L^{i_1}_x$ carries an invariant flat connection.}

\section{Smooth dynamical foliations from Ruelle}

This section contains results from Ruelle \cite{ruelle.diff.ergodic}, in their general form in subsection \ref{sec.ruelle.submanifolds}, and in the special case of foliations on compact $C^0$ manifolds in subsection \ref{sec.filtered.foliations}.  Subsection \ref{sec.tensor.stability} contains some key results restricting the possible values of $F^t_x$-invariant tensors.  The formulations of Ruelle's results are not exactly the same as in his paper.  In particular, we work with a flow, and we keep track of invariance by the full centralizer of this flow.  Proposition \ref{prop.perturbation.spectrum} in section \ref{sec.tensor.stability} requires a simple modification of the innards of Ruelle's proof.  The present section thus contains references to and some recapitulation of Ruelle, as well as some further proofs.

In this section, $( \{ \varphi^t \}, M, \mu)$ and the cocycle $\{F^t_x\}$ are as in the introduction, but the action is not assumed to be a measurable contraction; that is, the Lyapunov spectrum of $T$ has the general form
$$- \infty < \lambda^{(1)} < \cdots < \lambda^{(s)} < 0 \leq \lambda^{(s+1)} < \cdots < \lambda^{(m)}$$ 

Let $M_T \subset M$ be the $\mu$-conull set where the conclusions of the MET hold for the cocycle $T$.  
For $x \in M_T$, denote the filtration corresponding to the negative portion of the Lyapunov decomposition of $T_{\bf 0} \mathcal{E}_x \cong \BR^n$ by
$${\bf 0} \subset V_x^1 \subset \cdots \subset V_x^s \subseteq \BR^n. $$

A key feature of Ruelle's theorem is that it applies to a noninvertible system.  We are interested here in cocycles over invertible transformations or flows, and invertibility will be an essential assumption in section \ref{sec.tensor.stability}.  A crucial role will also be played, however, by Ruelle's perturbation theorem \cite[Thm 4.1]{ruelle.diff.ergodic}, which applies to arbitrary sequences of linear transformations that are sufficiently close to $\{ T^k_x \}$; see proposition \ref{prop.perturbation.spectrum} below.

%%  Denote the filtration corresponding to the Lyapunov decomposition of $T_{\bf 0} \mathcal{E}_x \cong \BR^n$ by
%% $${\bf 0} \subset V_x^1 \subset \cdots \subset V_x^s \subset \cdots \subset V^m_x = \BR^n $$
%% defined for almost-every $x \in M$.  

\subsection{Submanifolds in fibers}
\label{sec.ruelle.submanifolds}

Following are results on the existence of the Lyapunov-Ruelle submanifolds in the fibers of $\mathcal{E}$ and their invariance by the centralizer of $\{ F^t_x \}$.  
%The cocycle $\{F^t_x \}$ and $M_T \subseteq M$ are as above.

\begin{theorem}[Ruelle 1979]
\label{thm.ruelle.submanifolds}
Assume $\sup_{0 \leq t \leq 1} \ln^+ \Vert F^{\pm t}_x \Vert_q \in L^1(M,\mu)$.  Let $M_F \subseteq M$ be the set where the conclusion of the ergodic theorem holds for this function under the flow $\{ \varphi^t \}$.  Let $M_0 = M_F \cap M_T$.

For all $x \in M_0$, there is a nested family of $C^q$ injectively immersed submanifolds 
$${\bf 0} \subset \nu_x^{1} \subset \cdots \subset \nu_x^s 
%\subseteq \overline{B(\alpha(x))} 
\subset \mathcal{E}_x.$$
For $i = 1, \ldots, s,$ the submanifolds $\nu^i_x$ are tangent at ${\bf 0}$ to $V_x^i$, and are characterized by
$$\nu_x^i = \{ u \in \mathcal{E}_x \ : \ \limsup_{t \rightarrow \infty} \frac{1}{t} \ln {\Vert F^t_x u \Vert} \leq \lambda^{(i)} \}$$ 
%% They have the further properties 
%% \begin{enumerate}
%% \item $T_{\bf 0} (\nu_x^r) = V_x^r$
%% \item For $u,v \in \nu_x^r$, 
%% $$\lim \sup_t \frac{1}{t} \ln {\Vert F^t_x u - F^t_x v \Vert} \leq \lambda^{(r)}$$
%% \end{enumerate}

%% The germs of $\{ \nu_x^r \subset \mathcal{E}_x \}_{x \in M}$ are $\{ F^t \}$-invariant---that is, for each $t$, the intersection of $F^t_x(\nu_x^r)$ with $\nu^r_{\varphi^tx}$ is the closure in each of an open submanifold through ${\bf 0}$.

The union $\cup_x \nu^i_x$ is $F^t_x$-invariant: $F^t_x(\nu^i_x) = \nu^i_{\varphi^t(x)}$. 
\end{theorem}

Ruelle's original statement is for sequences rather than flows.  The following proposition explicates the relation between his theorem and the above statement modulo this difference.

\begin{proposition}
\label{prop.submanifolds.sequence}
The conclusions of the theorem hold for the sequences $\{F^k_x = F(k,x) \}$ and $\{ T^k_x = T(k,x) \}$ corresponding to iteration of $\varphi = \varphi^1$.  
\end{proposition}

\begin{Pf}
Let $\lambda^{(i)} < \lambda < \lambda^{(i+1)}$; if $i=s$, replace $\lambda^{(i+1)}$ with $0$.  Ruelle's theorem 5.1 a) \cite{ruelle.diff.ergodic} asserts the existence of measurable functions $\beta(x) > \alpha(x) > 0$ such that the sets 
$$\accentset{\circ}{\nu}_x^\lambda = \{ u \in B(\alpha(x)) \ : \ \Vert F^k_x (u) \Vert \leq \beta(x) e^{k \lambda} \ \forall k \geq 0 \}$$
are $C^q$ submanifolds of $\mathcal{E}_x$, tangent at ${\bf 0}$ to $V^i_x$, for all $x$ in a $\mu$-conull subset $\Gamma$.  The set $\Gamma$ is given by \cite{ruelle.diff.ergodic} (5.2)--(5.4) and, under our hypotheses, contains $M_0$.  
%By point b') of that theorem, these submanifolds are independent of $\lambda \in (\lambda^{(i)}, \lambda^{(i+1)})$, up to germs around ${\bf 0}$ in $\mathcal{E}_x$.

Remark 5.2 c) of \cite{ruelle.diff.ergodic} (see also (6.2) on p 52) concerns an important lower bound on the shrinkage of $\alpha$ along $\{ \varphi^k \}$-orbits: given $\lambda^{(s)} < \zeta < 0$, one can arrange that $\alpha(\varphi^k(x))$ decreases less rapidly than $e^{k \zeta}$.  
%% Thus for each $u \in \nu_x^r$, there is $k_0$ for which $F^k_x(u) \in B(\alpha(\varphi^k(x))$, and therefore in $\nu_{\varphi^k(x)}^r$, for all $k \geq k_0$.  Therefore, if $\tau$ is a tensor on $\nu_x^r$, then we can take limits for $u \in \nu_x^r$ as $k \rightarrow \infty$ involving $(F^k_x)_* \tau(u)$. 

Set 
$$\nu^i_x = \bigcup_{k \geq 0} (F^k_x)^{-1}(\accentset{\circ}{\nu}_{\varphi^k(x)}^\lambda)$$
For $u \in \intnu_x^\lambda$ and $\lambda^{(i)} < \lambda'< \lambda$, part b') of \cite[Thm 5.1]{ruelle.diff.ergodic} gives 
$$\limsup_{k \rightarrow \infty} \frac{1}{k} \ln \Vert F^k_x u \Vert \leq \lambda'$$
  It follows that
$$ \nu_x^i \subseteq \{ u \in \mathcal{E}_x \ : \ \limsup_{k \rightarrow \infty} \frac{1}{k} \ln {\Vert F^k_x u \Vert} \leq \lambda^{(i)} \}$$
Proof of the reverse containment follows an argument on the bottom of p 52 in \cite{ruelle.diff.ergodic}: if $u$ belongs to the right-hand set above, then for $\lambda^{(i)} < \lambda < \min\{ \lambda^{(i+1)},\zeta \}$, and sufficiently large $C$, $\Vert F^k_x (u) \Vert \leq C e^{k \lambda}.$ 
For sufficiently large $l$, $\alpha(\varphi^l(x)) > C e^{l \zeta}$.  Then, for all $k \geq 0$,
$$ \Vert F^{k+l}_x (u) \Vert \leq C e^{l \zeta} e^{k \lambda} < \alpha(\varphi^l(x)) e^{k \lambda} < \beta(\varphi^l(x)) e^{k \lambda}$$
so $F^l_x(u) \in \accentset{\circ}{\nu}_{\varphi^l(x)}^\lambda$.  Therefore $u \in \nu_x^i$.

Now $F^k_x$-invariance (both forwards and backwards) of the $\cup_x \nu^i_x$ follows easily from the $\lim \sup$ characterization. 

%% , which is set in a $C^q$ manifold where $\varphi$ is a diffeomorphism and $T^k_x$ is the derivative cocycle, restricted to the sum of negative Lyapunov spaces.  The required argument, which can be found , does not use any hypotheses particular to that setting, and is equally valid in our setting. 

%% The $F^k_x$-invariance of germs is easy.  It is clear that 
%% $$ F^k_x(\nu_x^i) \subset \{ u \in \mathcal{E}_{\varphi^k(x)} \ : \ \lim \sup_m \frac{1}{m} \ln {\Vert F^m_{\varphi^k(x)} u \Vert} \leq \lambda^{(i)} \}.$$

%% The image $B =F^k_x(\overline{B(\alpha(x))})$ contains ${\bf 0}$ and is diffeomorphic to a closed ball in $\mathcal{E}_{\varphi^k(x)}$. The intersection 
%% $$ F^k_x(\nu_x^i) \cap \nu_{\varphi^k(x)}^i = F^k_x(\nu^i_x) \cap \overline{B(\alpha(\varphi^k(x))}) = B \cap \nu_{\varphi^k(x)}^i$$  

On p 53 of \cite{ruelle.diff.ergodic} is the further argument showing $\nu_x^i$ is an injectively immersed submanifold: namely, each $(F^m_x)^{-1}(\accentset{\circ}{\nu}_{\varphi^m(x)}^\lambda)$ is contained in $(F^{l_0}_x)^{-1}(\accentset{\circ}{\nu}_{\varphi^{l_0}(x)}^\lambda)$ for some $l_0 > m$, in fact, for all $l \geq l_0$.  Now $\nu_x^i$ is an increasing union of $C^q$ submanifolds which are topological disks, tangent at ${\bf 0}$ to $V_x^i$; the remaining conclusions follow.
\end{Pf}

\begin{remark}
\label{rmk.use.ruelle}
In several key places below, reference will be made to the smaller submanifolds $\accentset{\circ}{\nu}_x^\lambda$ in the foregoing proof, and the functions $\alpha$ and $\beta$ in their definitions.  In a few places, we will have to recall and tinker with some details of their construction.
\end{remark}

\begin{remark} 
\label{rmk.jets.arb.small}
By part b) of \cite[Thm 5.1]{ruelle.diff.ergodic}, the derivatives of $F^k_x$ become arbitrarily small for sufficiently large $k$, uniformly on compact subsets of $\accentset{\circ}{\nu}_x^s$.  For $i \leq q -1$, consider the cocycle $\tilde{F}^{(i)}(k,x) \in \mbox{Diff}(\BR^n \times \BR^n, {\bf 0})$ mapping $(u,v)$ to $(F^k_x(u), \mathcal{S}(D^{(i)}_u F^k_x)(v))$, where $\mathcal{S}$ converts an $i$-symmetric tensor to the corresponding homogeneous polynomial of degree $i$.  Then $\tilde{F}^{(i)}$ satisfies the finite first moment condition of theorem \ref{thm.ruelle.submanifolds} with respect to the $C^{q-i}$-norm.  (See section \ref{subsec.Tr.spectrum} for a similar but somewhat more demanding verification.)  Then \cite[Thm 5.1 b)]{ruelle.diff.ergodic} implies that $D_u (\mathcal{S} (D^{(i)} F^k_x))$ becomes arbitrarily small uniformly in $u \in \accentset{\circ}{\nu}_x^s$, and it follows that $D^{(i+1)} F^k_x$ grows arbitrarily small for sufficiently large $k$, uniformly on compact subsets of $\accentset{\circ}{\nu}_x^s$.  
\end{remark}

The following proposition asserts that a bundle automorphism centralizing the initial cocycle preserves its Lyapunov-Ruelle submanifolds.  We first define the centralizer.

\begin{definition}
\label{def.centralizer}
Let $\{ F^k_x \}$ satisfy the discrete-time version of the standing assumptions of this section.
%, as well as $\ln^+ \Vert F^1_x \Vert_q \in L^1(M,\mu)$.  
The \emph{centralizer of $\{ F^k_x \}$}, denoted $Z(F^k_x)$, comprises all $(\{ G^k_x \},\psi)$, where 
\begin{enumerate} 
\item $\psi$ is a measurable automorphism of $(M,\mu)$ commuting with $\varphi$
\item $G(k,x) = G^k_x \in \mbox{Diff}^q(\BR^n,{\bf 0})$ is a measurable cocycle over $\psi$ satisfying for all $k,l \in \BZ$,
$$ G(k,\varphi^l(x)) \circ F(l,x) = F(l,\psi^k(x)) \circ G(k,x)$$
\item $\ln^+ \Vert G^{\pm 1}_x \Vert_q \in L^1(M,\mu)$.
\end{enumerate} 
The \emph{centralizer of $\{ F^t_x \}$}, $t \in \BR$, is defined similarly, and denoted $Z(F^t_x)$.
\end{definition}

\begin{proposition}
\label{prop.submanifolds.centralizer}
Let $\varphi$, $M_0$, and $\{ \nu_x^i \}$ be as in proposition \ref{prop.submanifolds.sequence}, and let $(\{ G^k_x \}, \psi) \in Z(F^k_x)$.  Let $M_G \subseteq M$ comprise the points for which the conclusion of the ergodic theorem holds for both functions $ \ln^+ \Vert G^{\pm 1}_x \Vert_q$ under the flow $\{ \varphi^t \}$.  Then $G^k_x$ preserves $\cup_x \nu^i_x$ for all $i = 1, \ldots, s$, for all $x \in \cap_{i \in \BZ} \psi^i(M_0 \cap M_G)$. 
%(** \emph{either intersect over all $\psi$ translates or verify that can do some kind of joint Birkhoff with commuting flows?} **)
%---that is, $\forall x \in M_F$, the intersection $G^k_x(\nu^r_x) \cap \nu^r_{\psi^k(x)}$ is the closure of an open submanifold through ${\bf 0}$ of each.
\end{proposition}

\begin{Pf}
%% The centralizing assumption gives
%%   $$F^k_{\psi(x)} = G^1_{\varphi^k(x)} F^k_x G^{-1}_{\psi(x)}.$$
  For $x \in M_0$, and $1 \leq i \leq s$, take $u \in \nu_x^i$, and let $u' = G^1_x(u)$.  By the centralizing assumption,
  \begin{eqnarray}
   \limsup_{k \rightarrow \infty} \frac{1}{k} \ln \Vert F^k_{\psi(x)}(u') \Vert  & = &  \limsup_{k \rightarrow \infty} \frac{1}{k} \ln \Vert G^1_{\varphi^k(x)} F^k_x(u)  \Vert
\label{cent.exp1}
    \end{eqnarray}

  %% As $G(1,x)$ is a diffeomorphism with inverse $G(-1,\psi(x))$, we can bound $\Vert G(1,x)u - G(1,x)v \Vert^{-1}$, for $u,v \in \overline{B(\alpha(x))}$, by a constant times $\Vert u - v \Vert^{-1}$.  Therefore (\ref{cent.exp2}) is bounded by
  %% \begin{eqnarray}
  %%   \label{cent.exp3}
  %%                  \lim \sup \frac{1}{n} \ln \frac{\Vert G(1,\varphi^n(x)) F(n,x)u - G(1,\varphi^n(x)) F(n,x)v \Vert}{\Vert u - v \Vert}
  %%   \end{eqnarray}
  
  From proposition \ref{prop.submanifolds.sequence},
  $$ \limsup_{k \rightarrow \infty} \frac{1}{k} \ln \Vert F^k_x(u) \Vert \leq  \lambda^{(i)};$$
  in particular, for sufficiently large $k$, the image $F^k_x(u) \in \overline{B(1)}$.  For such $k$, the norm $\Vert G^1_{\varphi^k(x)} F^k_x(u) \Vert$ is bounded by the $C^1$-norm $\Vert G^1_{\varphi^k(x)} \Vert_1$, times $\Vert F^k_x(u) \Vert$.  Then we can also bound 
  $$ \Vert G^1_{\varphi^k(x)} F^k_x(u) \Vert \leq  \Vert G^1_{\varphi^k(x)} \Vert_q \cdot \Vert F^k_x(u) \Vert$$

  Now (\ref{cent.exp1}) is bounded above by
  \begin{eqnarray}
    \label{cent.exp4}
   \limsup_{k \rightarrow \infty} \frac{1}{k} \ln  \left( \Vert G^1_{\varphi^k(x)} \Vert_q \cdot \Vert F^k_x(u) \Vert \right)
  \end{eqnarray}
  For $x \in M_G$,
  %% $$ \lim_k \frac{1}{k} \sum_{m=0}^{k-1} \ln^+ \Vert G^1_{\varphi^m(x)} \Vert_q = \int \ln^+ \Vert G^1_x \Vert_q\ d\mu(x)$$
 %% from which
    $$ \limsup_{k \rightarrow \infty}  \frac{1}{k} \ln \Vert G^1_{\varphi^k(x)} \Vert_q \leq \lim_{k \rightarrow \infty} \frac{1}{k} \ln^+ \Vert G^1_{\varphi^k(x)} \Vert_q = 0$$
  Therefore (\ref{cent.exp4}) is bounded by $\lambda^{(i)}$,
%   $$ \limsup_k \frac{1}{k} \ln \Vert F^k_x(u)  \Vert \leq  \lambda^{(i)}$$
as is (\ref{cent.exp1}), 
%is bounded by $\lambda^{(i)}$
  %% $$ \lim \sup \frac{1}{k} \ln \Vert F^k_{\psi(x)}(u') - F^k_{\psi(x)}(v') \Vert \leq \lambda^{(r)}$$
for any $u' \in G^1_x (\nu_x^i)$.  Thus $G^1_x(\nu_x^i) \subseteq \nu_{\psi(x)}^i$ for $x \in M_0 \cap M_G \cap \psi^{-1} (M_0)$.  

A similar argument gives $G^{-1}_x(\nu_x^i) \subseteq \nu_{\psi^{-1}(x)}^i$ for $x \in M_0 \cap M_G \cap \psi(M_0)$.  Then for $x \in \cap_{i \in \BZ} \psi^i(M_0 \cap M_G)$, we conclude $G^k_x(\nu_x^i) = \nu_{\psi(x)}^i$ for all $k \in \BZ$. \end{Pf}
%\{ u \in \mathcal{E}_{\psi(x)} \ : \ \limsup_k \frac{1}{k} \Vert F^k_{\psi(x)} (u) \Vert \leq \lambda^{(r)} \}.$$
%% As in the proof of lemma \ref{lemma.submanifolds.sequence}, 
%% $$ G^1_x(\nu_x^r) \cap \nu_{\psi(x)}^r = G^1_x(\nu_x^r) \cap \overline{B(\alpha(\psi(x)))} = \nu_{\psi(x)}^r \cap G^1_x(\overline{B(\alpha(x))}$$

Now we can finish the

\begin{Pf}(of Theorem \ref{thm.ruelle.submanifolds})
Take the submanifolds $\nu_x^i$ associated to $\{ F^k_x \}$, as in proposition \ref{prop.submanifolds.sequence}.  The $F^t_x$-invariance for all $t \in \BR$ follows from proposition \ref{prop.submanifolds.centralizer}; note that $M_0$ is $\varphi^t$-invariant.

To complete the proof, it suffices to show that for $u \in \nu_x^i$, 
$$ \limsup_{t \rightarrow \infty} \frac{1}{t} \ln \Vert F^t_x(u) \Vert \leq \lambda^{(i)}$$

Any $t$ can be written $n_t + \epsilon_t$ with $n_t \in \BN$ and $0 \leq \epsilon_t < 1$.  For sufficiently large $t$, the image $F^{n_t}_x(u) \in \overline{B(1)}$. Then 
$$\Vert F^t_x(u) \Vert  =  \Vert F^{\epsilon_t}_{\varphi^{n_t}(x)} F^{n_t}_x (u) \Vert \leq \Vert F^{\epsilon_t}_{\varphi^{n_t}(x)} \Vert_q \cdot \Vert F^{n_t}_x(u) \Vert$$
%is bounded by the $C^1$ norm of $F^{\epsilon_t}_{\varphi^{n_t}(x)}$ on $\overline{B(1)}$, times $\Vert F^{n_t}_x(u) - F^{n_t}_x(v) \Vert$, so it is also bounded by.  
Now
$$ \frac{1}{t} \ln \Vert F^t_x(u) \Vert \leq \frac{1}{n_t}  \sup_{0 \leq \epsilon \leq 1}  \ln \left( \Vert F^{\epsilon}_{\varphi^{n_t}(x)} \Vert_q \cdot \Vert F^{n_t}_x(u) \Vert \right)$$

As in the centralizer proposition, the ergodic theorem gives 
$$\lim_{k \rightarrow \infty} \frac{1}{k} \sup_{0 \leq \epsilon \leq 1} \ln^+ \Vert  F^{\epsilon}_{\varphi^{k}(x)} \Vert_q = 0$$
so 
$$ \limsup_{t \rightarrow \infty} \frac{1}{t} \ln \Vert F^t_x(u)  \Vert \leq \limsup_{t \rightarrow \infty} \frac{1}{n_t} \ln \Vert F^{n_t}_x(u) \Vert \leq \lambda^{(i)}$$
\end{Pf}

\begin{remark}
\label{remark.cinfinity.can}
On regularity: If $F^t_x$ acts by $C^\infty$ diffeomorphisms of $\BR^n$ and satisfies the finiteness condition that
$$ \sup_{0 \leq t \leq 1} \ln^+ \Vert F^{\pm t}_x \Vert_q \in L^1(M, \mu) \qquad \mbox{for all} \ q \geq 1,$$
then the submanifolds $\nu^r_x$ in theorem \ref{thm.ruelle.submanifolds} are $C^\infty$.  Section 5.3 of \cite{ruelle.diff.ergodic} verifies the smoothness of the $\accentset{\circ}{\nu}_x^\lambda$. (Note that the conull set $M_0$ can vary between the $C^q$ and $C^\infty$ cases).  The arguments in the proof of proposition \ref{prop.submanifolds.sequence} are also valid in the $C^\infty$ setting.  The hypothesis of proposition \ref{prop.submanifolds.centralizer} should then be replaced with $G^k_x \in \mbox{Diff}^\infty(\BR^n,{\bf 0})$, and $\ln^+ \Vert G^{\pm 1}_x \Vert_q \in L^1(M,\mu)$ for all $q \geq 1$.  We note that there is also a $C^\omega$ version of these submanifolds (see \cite[Sec. 5.4]{ruelle.diff.ergodic}).

Ruelle's regularity assumption is actually that the $F^k_x$ are of class $C^{q,\theta}$ with $q \geq 1$ and $\theta \in (0, 1]$---that is, the derivatives of order $q$ are H\"older continuous of exponent $\theta$.  The submanifolds $\nu_x^i$ are then also $C^{q,\theta}$-regular.  The results of this section can be verified in this greater generality.  See also \cite{kalinin.sadovskaya.nonuniform.normal}.  We will, however, use $q \geq 2$ for theorems \ref{thm.main} and \ref{thm.foliation.atlas} below.  
%We chose to avoid the additional notational burden. 
\end{remark}

\begin{remark}
One need not assume that $\mu$ is $\varphi^t$-ergodic.  An arbitrary invariant probability measure can be decomposed into ergodic components.  We make this assumption mostly to simplify statements and proofs.
\end{remark}

\subsection{Foliations by filtered stable manifolds}
\label{sec.filtered.foliations}

We now let $M$ be a compact $C^0$ manifold and $\{ \varphi^t \}$ a $C^0$ flow.  We will assume that the fibers of $\mathcal{E}$ are plaques of a $C^0$ foliation $L$ of $M$ by $n$-dimensional submanifolds, admitting a $\varphi^t$-invariant $C^q$ smooth structure.  Thus $M$ admits a $C^0$ foliated atlas in which the transitions are also $C^q$ along the leaves of $L$.  Fix a $C^0$ norm on $TL \rightarrow M$; we will assume $\{ \varphi^t \}$ is nonuniformly contracting along $L$.

We will call a \emph{$C^q$ atlas along $L$} a family of $C^q$ smooth charts $\theta_x : (\mathcal{E}_x,{\bf 0}) \rightarrow (L_x,x)$, varying measurably in $x$, and with $C^q$ norms of $\theta_x$ on $\overline{B(1)}$ and of $(\left. \theta_x \right|_{\overline{B(1)}})^{-1}$ bounded in $x$.
The cocycle $\{ F^t_x \}$ is the local action of $\{ \varphi^t \}$ on $L_x$ in these charts: 
$$ F^t_x = (\theta_{\varphi^t(x)})^{-1} \circ \varphi^t \circ \theta_x.$$  
Each $F^t_x$ is defined on a neighborhood of the origin in $\mathcal{E}_x \cong \BR^n$.

A $C^q$ atlas along $L$ will be called \emph{uniformly biLipschitz} if there is $\kappa \geq 1$ such that for some continuous metric $d$ on $M$, for almost all $x$,
$$ \frac{1}{\kappa} \Vert u - v \Vert \leq d(\theta_x(u), \theta_x(v)) \leq \kappa \Vert u - v \Vert \qquad \forall u,v \in \overline{B(1)}$$
It is not hard to find a uniformly biLipschitz atlas along $L$: take a finite collection of $C^0$ foliated charts $\varphi_i: V_i \rightarrow U_i, i = 1, \ldots, N,$ for which $V_i$ are open balls of radius $2$ and $\cup_i \varphi_i(\frac{1}{2} \overline{V}_i)$ form a finite cover of $M$.  Then make a piecewise continuous assignment $x \mapsto i$ with $x \in \varphi_i(\frac{1}{2} \overline{V}_i)$ and set $\theta_x = \varphi_i$ restricted to the leaf in $V_i$ through $\varphi^{-1}_i(x)$ and recentered at $x$.

%% An assumption of finite first moment is required:  
%% $$ \sup_{0 \leq t \leq 1} \ln^+ \Vert \left. \varphi^{\pm t} \right|_{L} \Vert_{q,x} \in L^1(M,\mu)$$ 
%% where $\Vert \cdot \Vert_{q,x}$ means the $C^q$ norm over the unit ball centered at $x$ (say, with respect to $d$).  
 Here is the main result of this section, which is very close to \cite[Thm 6.3]{ruelle.diff.ergodic}.
%% The goal of this subsection is to show that the leaves $L_x$ are foliated by a filtered family $L^1_x \subset \cdots L^s_x \subset W_x$ corresponding under $\theta_x$ to the submanifolds $\nu_x^1 \subset \cdots \subset \nu_x^s \subseteq \mathcal{E}_x$, and that these foliations are transversely $C^{q-1}$ inside each $L_x^s$.

\begin{theorem}[Ruelle 1979]
\label{thm.ruelle.foliations}
Let $\{ \varphi^t \}$ be a $C^0$, measure-preserving flow on a compact manifold $(M, \mu)$, preserving a $C^0$ foliation $L$ and a $C^q$ structure on the leaves.  Assume $\sup_{-1 \leq t \leq 1} J_x^{(q)} ( \left. \varphi^t \right|_L )$ is bounded in $x$ and that $D(\left. \varphi^t \right|_L)$ has all negative Lyapunov exponents in a uniformly biLipschitz $C^q$ atlas $\{ (\mathcal{E}_x, \theta_x) \}$ along $L$.  Then the corresponding cocycle satisfies the hypotheses of theorem \ref{thm.ruelle.submanifolds}; the resulting submanifolds $\nu_x^1 \subset \cdots \subset \nu_x^s = \mathcal{E}_x$ satisfy, for $x \in M_0$:

\begin{enumerate}
\item  Suppose $L_x = L_y$ and $y \in M_0 \cap \theta_x(\nu_x^i)$ for some $1 \leq i \leq s$.  Then $\theta_y(\nu_y^i) \cap \theta_x(\nu_x^i)$ is relatively open in each term.
\item The resulting foliations $L^1 \subset \cdots \subset L^s = L$ are $\varphi^t$-invariant, and can be defined at all points of $L_x^s$, for all $x \in M_0$.
\item  The distributions $\mathcal{V}_x^i = T_x (L^i_x)$ vary $C^{q-1}$-smoothly within $L_x^s$ for $1 \leq i < s$.
\end{enumerate}
\end{theorem}

\begin{Pf}
The assumption of negative Lyapunov exponents along $L$ implies $\{ \varphi^t \}$ is nonuniformly contracting along $L$, so $\{ F^t_x \}$ is defined on $\overline{B(1)}$ for sufficiently large $t$. The boundedness assumption on $q$-jets of $\varphi^t$ along $L$ implies $\{ F^t_x \}$ satisfies the finite first moment assumption of theorem \ref{thm.ruelle.submanifolds}. 

Now recall that $\nu_x^s$ is the increasing union $\cup_{k \geq 0} (F^k_x)^{-1} \intnu_{\varphi^{k}(x)}^\lambda$ for $\lambda^{(s)} < \lambda < 0$.  For any compact $D \subset \nu_x^s$, there is $k_0$ such that $F^k_x$ is guaranteed to be defined on $D$ for all $k \geq k_0$. 

Let $\kappa$ be the biLipschitz constant relating $\{ \theta_x \}$ with a distance $d$ on $M$.  Let $u \in \nu_x^i$ with $\theta_y^{-1} \theta_x(u) = v$.  Then 
\begin{eqnarray*}
\Vert F^k_y(v) \Vert & = & \Vert \theta_{\varphi^k(y)}^{-1} \circ \varphi^k(\theta_x(u)) \Vert \\
& = & \Vert \theta_{\varphi^k(y)}^{-1} \circ \theta_{\varphi^k(x)} \circ F^k_x (u) \Vert \\
& \leq & \kappa  \left( d(\varphi^k(y),\varphi^k(x))+  d(\theta_{\varphi^k(x)} \circ F^k_x (u),\varphi^k(x)) \right) \\
& \leq & \kappa^2 \left( \Vert F^k_x(\theta_x^{-1}(y)) \Vert  + \Vert F^k_x(u) \Vert \right)
\end{eqnarray*}
so
$$
\limsup_{k \rightarrow \infty} \frac{1}{k} \ln \Vert F^k_y(v) \Vert  \leq \limsup_{k \rightarrow \infty} \frac{1}{k} \ln \kappa^2 \max\{ \Vert F^k_x(\theta_x^{-1}(y)) \Vert , \Vert F^k_x (u) \Vert \} \leq \lambda^{(i)}$$

From the characterization in theorem \ref{thm.ruelle.submanifolds}/proposition \ref{prop.submanifolds.sequence}, we conclude $v \in \nu_y^i$.  Thus
$$ \theta_x(\nu_x^i) \cap \theta_y(\mathcal{E}_y) = \theta_x(\nu_x^i) \cap \theta_y(\nu_y^i)$$
which, similarly, equals $\theta_y(\nu_y^i) \cap \theta_x(\mathcal{E}_x)$.
%$$ \theta_x(\nu_x^i) \cap \theta_y(\nu_y^i) = \theta_y(\nu_y^i) \cap \theta_x(\mathcal{E}_x)$$
Now (1) is proved.

The $\varphi^t$-invariance of the foliations $L^i$, $i = 1, \ldots, s$, follows immediately from their definition in terms of the atlas $(\mathcal{E}_x,\theta_x)$ and the $F^t_x$-invariant submanifolds $\nu_x^i$.  Given $x \in M_0$ and $y \in L_x^s$, the trajectories $\varphi^t(x)$ and $\varphi^t(y)$ come arbitrarily close as $t \rightarrow \infty$: one can choose a path between them in $L_x^s$ and cover this path with finitely many open sets of the form $\theta_{x_i}(\intnu_{x_i}^s)$.  Each such open interval is uniformly exponentially contracted under $\varphi^t$ as $t \rightarrow \infty$. Let $t_0$ be a time at which the images of $x$ and $y$ are contained in the image of $\theta_{\varphi^{t_0}(x)}$.  All foliations are defined on this set, so $L^1_{\varphi^{t_0}(y)} \subset \cdots \subset L^s_{\varphi^{t_0}(y)}$ can be defined on a neighborhood of $\varphi^{t_0}(y)$.  Then to prove (2), let $L^i_y = \varphi^{-t_0}(L^i_{\varphi^{t_0}(y)})$ for each $1 \leq i < s$.  

Given $x \in M_0$ and $\lambda^{(1)} < \lambda_1 < \lambda^{(2)} < \cdots < \lambda = \lambda_s < 0$, let $\alpha$ be the minimum over $i = 1, \ldots, s$ of the radii $\alpha(x)$ in the definition of $\intnu_x^{\lambda_i}$. 
%(see the the proof of proposition \ref{prop.submanifolds.sequence}).  
We will show that the distribution $\mathcal{V}^i_y$ varies smoothly over $y \in \theta_x(B(\alpha)) \cap L_x^s$, for any $1 \leq i < s$.  Let $u = \theta_x^{-1}(y)$. Note that
$$\theta_x^{-1}(L^i_y) \subseteq \{ v \in  \mathcal{E}_x : \ \limsup_{t \rightarrow \infty} \frac{1}{t}\ln \Vert F^t_x(v) - F^t_x(u) \Vert \leq \lambda^{(i)} \}$$
and $\theta_x^{-1}(L^i_y)$ contains a neighborhood of $u$ in the right-hand set.  Remark 5.2 (b) of \cite{ruelle.diff.ergodic} shows that the tangent space at $u$ of the latter set varies $C^{q-1}$-smoothly in $u \in \accentset{\circ}{\nu}_x^\lambda$.  Mapping forward by $\theta_x$ gives (3). 
\end{Pf}

\begin{remark}
\label{rmk.foliations.centralizer}
Let $\{ \psi^k \}$ be a continuous transformation commuting with $\{ \varphi^t \}$, preserving the foliation $L$, acting $C^q$ differentiably along leaves, with $J_x^{(q)} ( \left. \psi \right|_L )$ bounded in $x$.  It follows from proposition \ref{prop.submanifolds.centralizer} that the foliations $L^1 \subset \cdots \subset L^s \subset L$ are $\psi^k$-invariant.
%, provided the function $\ln^+ \Vert \left. \psi^{\pm 1} \right|_L \Vert_{q,x}$ is also $\mu$-integrable.  
We will denote the collection of such $\psi$ by $Z(\varphi^t)$.
\end{remark}

% \begin{remark}
% $C^\infty$ case, $C^\omega$ case.
% \end{remark}

%% \begin{corollary}
%% Let ...  Then, there is a nested family of foliations
%% $$ W^1_x \subset \cdots W^s_x \subset W_x$$ 
%% varying $C^{q-1}$-smoothly in each $W_x$ and continuously in $M$.  These foliations are $\varphi^t$-invariant and tangent to the negative Lyapunov filtration for $D \varphi^t$ along $W$; they are moreover characterized by
%% $$ \lim \sup$$
%% with respect to any continuous distance $d$ on $M$.
%% \end{corollary}

\subsection{Stability of invariant tensors}
\label{sec.tensor.stability}

The aim of this section is to establish restrictions on values of $F^t_x$-invariant tensors on the Lyapunov-Ruelle manifolds $\nu_x^i$, for $1 \leq i \leq s$.  Sequences of the form $\{ D_u F^k_x \}_{k \in \BN}$, for $u \in \nu_x^i$, have the same spectrum as $\{ T^k_x \}_{k \in \BN}$, by Ruelle's perturbation theorem.  We will show that in fact the push-forward by $\{ F^k_x \}$ of a tensor at $u \in \nu_x^i$ has the same asymptotic expansion rates as the push-forward at ${\bf 0}$ by $\{ T^k_x \}$.  

\subsubsection{Digression on Ruelle's proof: constants in the construction}  
\label{subsubsec.digression}
Here we must expose some internal constants in Ruelle's construction of the $\intnu_x^\lambda$ (see remark \ref{rmk.use.ruelle}), in order to introduce the proof of proposition \ref{prop.perturbation.spectrum}.  

Ruelle's perturbation theorem  \cite[Thm 4.1]{ruelle.diff.ergodic} says that given a linear sequence $\{ T_k \}$ satisfying 
\begin{itemize}
\item $\limsup_{k \rightarrow \infty} \ln \Vert T_k \Vert \leq 0$;
\item $\lim_{k \rightarrow \infty} (T^{k*} T^k)^{1/2k} = \Lambda$ where $T^k = T_k \cdots T_1$;
\end{itemize}
and given $\eta > 0$, there exists $\delta > 0$, such that for any linear sequence $\{ T_k'\}$ with 
$$ \Vert T - T' \Vert = \sup_k \Vert T_k - T_k' \Vert e^{3 k \eta} < \delta,$$
$\{ T^{'k} \}$ has a well-defined spectrum, equal to that of $\{T^k \}$.  For $\lambda^{(i)}$ in this spectrum, and $P^{(i)}$ the projection to the corresponding asymptotic eigenspace for $\{ T^{'k} \}$, equation (4.5) of \cite{ruelle.diff.ergodic} additionally says that given $\epsilon > 0$, there is $B_{\epsilon}' > 1$ such that $\Vert T^{'k} \circ P^{(i)} \Vert \leq B_\epsilon' \cdot e^{k (\lambda^{(i)} + \epsilon)}$.

This theorem is applied several times in the construction of $\accentset{\circ}{\nu}_x^\lambda$, with constants defined as follows.  First, $\lambda = \lambda^{(i)} + \epsilon$, with $\lambda < \lambda^{(i+1)}$ (if $i = s$, take $\lambda^{(i+1)} = 0$).  Then $\eta$ is such that $0 < 4 \eta \leq - \lambda$.  Now $\delta$ is given by the perturbation theorem applied to the sequence $\{ T_k = T(1,\varphi^{k-1}(x)) \}$ (He further requires $\delta < 1/\sqrt{2} A$, where $A$ is another constant given by \cite[Thm 4.1]{ruelle.diff.ergodic}).  Next, $0 < \beta(x) < \min\{ 1, \delta/G \}$, where
\begin{eqnarray}
\label{eqn.def.G}
 G(x) = \sup_k \Vert F^1_{\varphi^{k-1}(x)} \Vert_q \cdot e^{- k \eta - \lambda}
\end{eqnarray}
Note $G(\varphi^l(x)) \leq e^{l \eta} G(x)$.  Last, $\alpha(x) = \beta(x)/B_\epsilon'$.  End of digression.
\bigskip

Before proceeding, we state some basic facts about behavior of Lyapunov exponents under restriction to invariant subspaces and quotients.

\begin{proposition}
\label{prop.basic.lyapunov}
Let $\{ T^t_x \}$ be a linear cocycle over an ergodic, measure-preserving flow $( \{ \varphi^t \}, M, \mu)$, satisfying (\ref{eqn.met}).  Let $\Sigma$ be the Lyapunov spectrum of $\{ T^t_x \}$, with Lyapunov decomposition $W^1_x \oplus \cdots \oplus W^m_x$.  Suppose that $\{ U_x \}_{x \in M}$ form a measurable $T^t_x$-invariant subbundle.  Then
\begin{enumerate}
\item The spectrum of $\{ T^t_x \}$ restricted to $U$ is a subset of $\Sigma$, and $U_x = \oplus_i (U_x \cap W^i_x)$.  
\item The spectrum of $\{ \overline{T}^t_x \}$, the cocycle on the quotient by $U$, is a subset of $\Sigma$, and each $W^i$ maps surjectively to the corresponding Lyapunov distribution in the quotient.
\end{enumerate}
\end{proposition}

Item (1) follows quickly from the definition of the Lyapunov decomposition for an invertible system.  For (2), the Oseledec-Pesin reduction theorem \cite[Thm 6.10]{barreira.pesin.introset} is helpful.  It gives tempered equivalences (see \cite[p 103]{barreira.pesin.introset}) between $\{ T^1_x \}$ and block $\epsilon$-approximately conformal matrices.  The blocks correspond to the Lyapunov decomposition.  If $W^i_x \cap U_x = 0$, then the angle $\angle(W^i_{\varphi^t(x)}, U_{\varphi^t(x)})$ decreases subexponentially as $t \rightarrow \infty$.  Thus for $v \in W^i_x$, the projection of $T^t_x v$ to $U_{\varphi^t(x)}^\perp$ has norm shrinking like $e^{t \lambda^{(i)}}$.  For ${\bf 0} \neq W^i_x \cap U_x \neq W^i_x$, the restriction of $T^t_x$ to $W^i_x$ is tempered equivalent to an $\epsilon$-approximately scalar matrix, for any $\epsilon$, which again implies that $\angle(W^i_{\varphi^t(x)} \cap U_{\varphi^t(x)},W^i_{\varphi^t(x)} \cap T^t_x(U_x^\perp))$ decreases subexponentially. 

Now let $\mathcal{T} \rightarrow M$ be a tensor bundle derived from $T_{{\bf 0}}\mathcal{E} \rightarrow M$ via $T^t_x$-invariant subbundles, quotients, and tensor operations.  Write $\{ T^t_{x*} \}$ for the linear cocycle induced from $\{T^t_x \}$ on $\mathcal{T}$.  Note that as $\{ T^t_x \}$ satisfies (\ref{eqn.met}), so does $\{ T^t_{x*} \}$.  Let $\Sigma$ be the spectrum of $\{ T^t_{x*} \}$ on $\mathcal{T}$, and denote by $\mathcal{T}_x^{\sigma}$ the Lyapunov distribution corresponding to $\sigma \in \Sigma$.

\begin{proposition}
\label{prop.T.stability}
If $\tau : M \rightarrow \mathcal{T}$ is a measurable, $T^t_{x*}$-invariant section, then \ $\tau(x) \in \mathcal{T}^{0}_x$ for almost-every $x \in M_0$.
\end{proposition}

\begin{Pf}
Suppose that $\tau(x)$ has a nontrivial component on $\mathcal{T}^\sigma_x$ for some $\sigma > 0$.  Then $\Vert (T^t_x)_* (\tau(x)) \Vert  = \Vert \tau(\varphi^t(x)) \Vert \rightarrow \infty$ as $t \rightarrow \infty$.  On the other hand, $\Vert \tau \Vert$ agrees with a continuous and bounded function on a set of measure $1-\epsilon$ by Lusin's theorem.  By ergodicity, for almost every $x$, the trajectory $\varphi^t(x)$ visits this set infinitely many times.  We conclude that for almost every $x$, the component of $\tau(x)$ on $\mathcal{T}^\sigma_x$ is zero.

If $\tau(x)$ has a nontrivial component on $\mathcal{T}^\sigma_x$ with $\sigma < 0$, then the same argument with $t \rightarrow - \infty$ leads to the same contradiction.
\end{Pf}

\begin{proposition}
\label{prop.perturbation.spectrum}
%Let $F(t,x)$ be as in theorem \ref{ruelle.submanifolds}, and $T(t,x) = D_{\bf 0} F^t_x$.  Let $\nu_x^r$ be given by theorem \ref{ruelle.submanifolds}, for some $1 \leq r \leq s$.  
Let $\mathcal{S} \rightarrow \cup_{x \in M_0} \nu_x^i $ be a tensor bundle constructed from $\oplus_x T \nu_x^i$ via $F^k_x$-invariant subbundles, quotients, and tensor operations, smooth in each $\nu^i_x$ and measurable in $x$.  Denote $\{ T_{\mathcal{S}}(k,x) \}$ the linear cocycle for $\{ T^k_{x*} \}$ on $\mathcal{S}_{\bf 0}$.  Fix $u \in \nu_x^i$, and consider the following sequence mapping $\mathcal{S}_u$ to $\mathcal{S}_{F^k_x(u)}$:
$$T'_{\mathcal{S}}(k) = (F^k_x)_{*u} = (F^1_{\varphi^{k-1}(x)})_{*F^{k-1}_x(u)} \cdots (F^1_x)_{*u}$$ 
Then $\{ T'_{\mathcal{S}}(k) \}$ has a well-defined spectrum, equal to that of $\{T_{\mathcal{S}}(k,x) \}$. 
\end{proposition}

\begin{remark}
In the proof below, we will modify $G$ in Ruelle's proof (see \cite[(5.5)]{ruelle.diff.ergodic}).  The functions $\beta(x)$ and $\alpha(x)$ will be modified accordingly; see subsection \ref{subsubsec.digression}.  The bound on shrinkage along orbits $e^{k \zeta} = O(\alpha(\varphi^k(x))$, for $\lambda^{(s)} < \zeta < 0$, will remain intact, and the $\{ \nu_x^i \}$ will be unchanged.  
\end{remark}

\begin{Pf}
It suffices to prove the statement when $u \in \accentset{\circ}{\nu}_x^\lambda$, for $\lambda$ as in subsection \ref{subsubsec.digression} above, because the spectrum of $T'_{\mathcal{S}}(k)$ is not changed by precomposition with an invertible linear map.
By definition, $\Vert F^k_x(u) \Vert \leq \beta(x) e^{\lambda k}$ for all $k \geq 0$.  As $q \geq 2$,
$$ \Vert D_{F^{k-1}_x(u)} F^1_{\varphi^{k-1}(x)} - T^1_{\varphi^{k-1}(x)} \Vert \leq \Vert F^1_{\varphi^{k-1}(x)} \Vert_q  \beta(x) e^{(k-1) \lambda} $$
%\leq \sup_{0 \leq t \leq 1} \Vert F^t_x \Vert_q \beta(x) e^{k \lambda}$$
%for $k$ sufficiently large that $F^{k-1}_x(u) \in \overline{B(1)}$.  FOLLOWS FROM BETA(X) < 1
Similarly,
$$ \Vert (D_{F^{k-1}_x(u)} F^1_{\varphi^{k-1}(x)})^{-1} - (T^1_{\varphi^{k-1}(x)})^{-1} \Vert \leq \Vert (F^1_{\varphi^{k-1}(x)})^{-1} \Vert_q  \beta(x) e^{(k-1) \lambda} $$
Because $x \in M_F$,
$$ \lim_{k \rightarrow \infty} \frac{1}{k} \ln^+ \Vert (F^1_{\varphi^{k-1}(x)})^{\pm 1} \Vert_q = 0$$

Let $\mathcal{S}$ be a tensor bundle of type $(a,b)$.
Let $\epsilon$ and $\eta$ be as in subsection \ref{subsubsec.digression}.  Let $\delta$ be given by the perturbation theorem applied to the linear sequence $\{ T_k = T(1,{\varphi^{k-1}(x)})_* \}$, acting on $\mathcal{S}_{\bf 0}$.  Instead of the function in (\ref{eqn.def.G}), set
$$G(x) = \sup_k \ (\Vert F^1_{\varphi^{k-1}(x)} \Vert_q)^b \cdot (\Vert (F^1_{\varphi^{k-1}(x)})^{-1} \Vert_q)^a \cdot e^{-k \eta - \lambda}  < \infty$$

Then choose $\beta(x) < \min\{ \delta/((a+b) G(x)), 1 \}$.  Now
\begin{equation*}
\begin{split}
&  \Vert (F^1_{\varphi^{k-1}(x)})_{*F^{k-1}_x(u)} - (T^1_{\varphi^{k-1}(x)})_* \Vert \\
\leq  b & \cdot \Vert D_{F^{k-1}_x(u)} F^1_{\varphi^{k-1}(x)} - T^1_{\varphi^{k-1}(x)} \Vert \cdot \Vert F^1_{\varphi^{k-1}(x)} \Vert_q^{b-1} \cdot \Vert (F^1_{\varphi^{k-1}(x)})^{-1} \Vert^a_q, \\
& + a \cdot \Vert (D_{F^{k-1}_x(u)} F^1_{\varphi^{k-1}(x)})^{-1} - (T^1_{\varphi^{k-1}(x)})^{-1} \Vert \cdot \Vert F^1_{\varphi^{k-1}(x)} \Vert_q^b \cdot \Vert (F^1_{\varphi^{k-1}(x)})^{-1} \Vert^{a-1}_q  \\
%% & \leq &   \Vert (D_{F^{k-1}_x(u)} F^1_{\varphi^{k-1}(x)})^{-1} - (T^1_{\varphi^{k-1}(x)})^{-1} \Vert^a \cdot \Vert D_{F^{k-1}_x(u)} F^1_{\varphi^{k-1}(x)} - T^1_{\varphi^{k-1}(x)} \Vert^b \\
& \leq    (a + b) \cdot \Vert (F^1_{\varphi^{k-1}(x)})^{-1} \Vert_q^a  \cdot \Vert F^1_{\varphi^{k-1}(x)} \Vert_q^b  \cdot \beta(x) e^{(k-1) \lambda} 
\end{split}
\end{equation*}
Then (see \cite[(5.5)]{ruelle.diff.ergodic}),
$$ \Vert (F^1_{\varphi^{k-1}(x)})_{*F^{k-1}_x(u)} - (T^1_{\varphi^{k-1}(x)})_* \Vert e^{3 k \eta} \leq (a + b) \cdot \beta(x) G < \delta $$
Now by \cite[Thm 4.1]{ruelle.diff.ergodic}, $\{ T_{\mathcal{S}}(k,x) \}$ and $\{ T'_{\mathcal{S}}(k) \}$ have the same spectrum.  As in (\ref{eqn.def.G}), $G(\varphi^l(x)) \leq e^{l \eta} G(x)$.  Then $\alpha$ and $\beta$ decrease along $\varphi$-orbits as in \cite[Rmk 5.2 c)]{ruelle.diff.ergodic} (see also \cite[Sec 4.7, p 43]{ruelle.diff.ergodic}).  Referring to the proof of proposition \ref{prop.submanifolds.sequence}, one can now see that our modifications of $\alpha$ and $\beta$ do not ultimately alter the manifolds $\{ \nu_x^i \}$.
\end{Pf}

\begin{proposition}
\label{prop.F.stability}
% Assume $\sup_{0 \leq t \leq 1} \ln^+ \Vert F^{\pm t}_x \Vert \in L^1(M,\mu)$.  
%Let $\nu_x = \nu_x^i$ be given by theorem \ref{thm.ruelle.submanifolds}, for some $1 \leq i \leq s$. 
Let $\oplus_x \tau_x$ be a family of smooth tensors belonging to a bundle $\mathcal{S} \rightarrow \cup_{x \in M_0} \nu_x^i$, varying measurably in $x$.  Assume that $\mathcal{S}$ is $F^k_x$-invariant, as in proposition \ref{prop.perturbation.spectrum}, and that $\tau$ is $F^t_x$-invariant---that is, 
$$(F^t_x)_* (\tau_x(u)) = \tau_{\varphi^t(x)}(F^t_x u).$$  
Let $\Sigma$ be the spectrum of $T_{\mathcal{S}}$ on $\mathcal{S}_{\bf 0}$, and denote $\Sigma^{\leq 0}$ the subset of nonpositive Lyapunov exponents. Then, for all $x$ in a $\varphi^t$-invariant, $\mu$-conull subset, 
%of $M_0$, 
$$ \lim_{t \rightarrow \infty} \frac{1}{t} \ln \Vert (F^t_x)_* (\tau_x(u)) \Vert \in \Sigma^{\leq 0} \qquad \forall \ u \in \nu_x^i.$$
\end{proposition}

\begin{Pf}
%% We will assume that the conclusion of the Birkhoff pointwise ergodic theorem holds for $\{ \varphi^t \}$ on $M_0$.  Consider the sequence of matrices, as on p 47 of \cite{ruelle.diff.ergodic}, $T_k' = DF_{\varphi^{k-1}(x)}(F^{k-1}_x(u))$.  Then (see p 48 of \cite{ruelle.diff.ergodic}), $|| T - T'||$ is less than the $\delta >0$ given by Ruelle's perturbation theorem \cite[Thm 4.1]{ruelle}.  By point (4.3) of this theorem, 
%% $$T^k = T_k \circ \cdots \circ T_1 = D_{\bf 0}F_{x}^k \ \mbox{and} \ (T')^k = T'_k \circ \cdots \circ T'_1 = D_uF_x^k$$
%% have the same spectrum; more precisely, 
%% $$ \lim_{k \rightarrow \infty}  (T^{'k*} T^{'k})^{1/2k} = \lim_{k \rightarrow \infty} (T^{k*}T^k)^{1/2k}$$ 
%% Then for any $v \in T_u \nu_x$, 
%% $$  \lim_{k \rightarrow \infty} \frac{1}{k} \ln \Vert D_u F^k_x v \Vert \in \{ \lambda^{(1)}, \ldots, \lambda^{(s)} \}$$
%% The action of $(F^k_x)_*$ on $\tau_x(u)$ is built from tensor powers of $D_u F^k_x$ and $(D_u F^k_x)^{-1}$, so has asymptotic growth rates given by integer combinations of $\lambda^{(1)}, \ldots, \lambda^{(s)}$.  These are the same as those for $T^k_x$, so they are in $\Sigma$.
By proposition \ref{prop.perturbation.spectrum} above, $\{ T'_{\mathcal{S}}(k) \}$ and $\{ T_{\mathcal{S}}(k,x) \}$ have the same spectrum, so
$$ \lim_{k \rightarrow \infty} \frac{1}{k} \ln \Vert (F^k_x)_{*u} (\tau_x(u)) \Vert \in \Sigma$$

To extend the limit to $t \rightarrow \infty$ in $\BR$, the argument is the same as in the proof of theorem \ref{thm.ruelle.submanifolds}.  
% Write $F^t_x = F^{\epsilon_t}_{\varphi^{n_t}(x)} \circ F^{n_t}_x$ with $0 \leq \epsilon_t < 1$, and compute that
% $$ \frac{1}{t} \ln \Vert (F^t_x)_{*u} (\tau_x(u)) \Vert \leq \frac{1}{n_t} \ln \Vert (F^{\epsilon_t}_{\varphi^{n_t}(x)})_{*F^{n_t}_x(u)} \Vert + \frac{1}{n_t} \ln \Vert (F^{n_t}_x)_{*u} (\tau_x(u)) \Vert  $$
% The first term vanishes in the limit as $t \rightarrow \infty$, and the second converges to a value in $\Sigma$.

For $\lambda^{(i)} < \lambda < \lambda^{(i+1)}$, consider the measurable function 
$$ S(x) = \sup_{u \in \accentset{\circ}{\nu}_x^\lambda} \Vert \tau_x(u) \Vert$$
 By Lusin's theorem, $S(x)$ agrees with a continuous and bounded function on a set of measure $1-\epsilon$, visited by $\varphi^t(x)$ for infinitely many $t$.  On the other hand, if 
$$  \lim_{t \rightarrow \infty} \frac{1}{t} \ln \Vert (F^t_x)_* (\tau_x(u)) \Vert > 0,$$
then $\Vert \tau_{\varphi^t(x)} (F^t_x u) \Vert \rightarrow \infty$ as $t \rightarrow \infty$.  This would be a contradiction.
\end{Pf}

\section{Geometric structures on fibers for contractions}
\label{sec.geom.str}

Here $(\{ \varphi^t \}, M,\mu)$ and $\mathcal{E} \rightarrow M$ are as in the previous section.  Also as above, $F(t,x) = F^t_x \in \mbox{Diff}^q(\BR^n,{\bf 0})$ and $T^t_x = D_{\bf 0} F(t,x)$ are the cocyles associated to the $\{ \varphi^t \}$-action on $\mathcal{E}$.  In this section the Lyapunov exponents $\Sigma^0$ of $\{T^t_x \}$ are assumed to \emph{all be negative}:
$$ - \infty < \lambda^{(1)} < \cdots \lambda^{(s)} < 0,$$
so $\{ F^t_x \}$ is asymptotically infinitesimally contracting on the fibers of $\mathcal{E}$ as $t \rightarrow \infty$.
We will use Ruelle's dynamical foliations to construct a family of $F^t_x$-invariant geometric structures on the fibers of $\mathcal{E}$.

\subsection{Prolongation of $F$}
\label{sec.prolongation}

Denote the general linear group of $\BR^n$ by $\GL(n)$.  The group $\GL^{(r)}(n)$ comprises the $r$-jets at ${\bf 0}$ of local diffeomorphisms of $\BR^n$ fixing ${\bf 0}$.  It can be identified with truncated Taylor series, or polynomials of degree $r$, with zero constant term and invertible linear component.  The Lie algebra $\alggl^{(r)}(n)$ comprises the $r$-jets at ${\bf 0}$ of local vector fields vanishing at ${\bf 0}$, and can be identified with the polynomials of degree $r$ with zero constant term.  The bracket of $X,Y \in \alggl^{(r)}(n)$ can be computed as $[X,Y] = J_{\bf 0}^{(r)} [\widetilde{X}, \widetilde{Y}]$, where $\widetilde{X}$ and $\widetilde{Y}$ are corresponding representative local vector fields.  We will denote $\rho^r_s$ the projection $\GL^{(r)}(n) \rightarrow \GL^{(s)}(n)$, for $r > s$, and also the projection $\alggl^{(r)}(n) \rightarrow \alggl^{(s)}(n)$.  We set $S^{(r)}(n) = \ker \rho^r_{r-1}$; it is isomorphic to the abelian group $\Sym^r(\BR^{n*}) \otimes \BR^n$.  

We will make use below of a norm on $\alggl^{(r)}(n)$.  View an element $X$ as an ordered $n$-tuple of polynomials of degree at most $r$, and define $\interleave X \interleave$ to be the maximum norm of the monomial coefficients in $X$.  Note that when $X$ is viewed as a smooth map of $\BR^n$, then this norm is bounded above by the $C^r$ norm $\Vert X \Vert_r$.  For $X \in \alggl(n)$, it is also bounded above by the usual linear operator norm (and below by $\Vert X \Vert /n$).
%, and write $X^{(l)}$ for the homogeneous component of degree $l$.  Define 
%$$ \Vert X^{(l)} \Vert = \sup_{\Vert v \Vert = 1} \Vert X^{(l)}(v) \Vert \qquad \mbox{and} \qquad \Vert X \Vert = \max_{1 \leq l \leq r} \Vert X^{(l)} \Vert$$  
We will need a bound on the growth of $\interleave \cdot \interleave$ under composition (followed by truncation above degree $r$).  
%% Write $P_r(n)$ for the sum of all monomials of degree at most $r$ in $n$ variables; it has $\binom{n+r}{r}$ terms.  
%% Let $c(r,n)$ be the maximum coefficient of a term in $P_r(n)^r$. GIVE A LITTLE EXPLANATION.

\begin{proposition}
\label{prop.poly.bounds}
Let $X, Y \in \alggl^{(r)}(n)$.  Then
$$ \interleave Y \circ X \interleave \leq \interleave Y \interleave \cdot \max \{ \interleave X \interleave, \interleave X \interleave^r \} \cdot c(r,n)$$
%% \item If $X$ is invertible, then 
%% $$ \Vert X \Vert \geq c(r,n)^{-1} \cdot \min \left\{ \frac{1}{\Vert X^{-1} \Vert}, \frac{1}{\Vert X^{-1} \Vert^r} \right\} $$
where $c(r,n)$ is a combinatorial constant depending on $r$ and $n$.
\end{proposition}

Observe that for a monomial $Y(u_1, \ldots, u_n)$ of degree $j$ and a polynomial $X = (X_1, \ldots, X_n)$, the maximum coefficient of $Y \circ X$ is bounded by $c'(j,n) \cdot \interleave Y \interleave \cdot \interleave X \interleave^j$, where $c'(j,n)$ is a combinatorial constant.  The claimed bound for arbitrary $Y,X \in \alggl^{(r)}(n)$ follows by a similar estimate.
%, and $1 \leq l \leq r$, the norm 
%$$ \Vert Y \circ X \Vert \leq \Vert Y \Vert \cdot \max  \{ \Vert X \Vert, \Vert X \Vert^l \} \cdot  c(r,n)$$
%for a combinatorial constant $c(r,n)$.
% depending on $l$ and $n$, increasing in $l$. Now taking the maximum over $1 \leq l \leq r$ gives the bound claimed in the proposition.

%% $Y = x_1^{j_1} \cdots x_n^{j_n}$ a monomial of degree $j_1 + \cdots + j_n = j$, 
%%  $$ \Vert Y \circ X \Vert \leq \max\{ \Vert X \Vert, \Vert X \Vert^r \} \cdot c(j,r,n) $$

In a manifold $M^n$, the order-$r$ frame bundle $\mathcal{F}^{(r)}(M)$ is a principal $\GL^{(r)}(n)$-bundle, comprising $r$-jets at ${\bf 0}$ of coordinate parametrizations of $M$.  The order-$r$ frame bundle of $\BR^n$ can be trivialized $\mathcal{F}^{(r)}(\BR^n) \cong \BR^n \times \GL^{(r)}(n)$ by identifying $(u,g)$ with the $r$-jet at ${\bf 0}$ of $g$ composed with translation by $u$.

Define the \emph{$r$th prolongation of $\mathcal{E}$} by $\mathcal{E}^{(r)}_x = \mathcal{F}^{(r)}(\mathcal{E}_x).$  Denote by $\pi^r_s$ the projection $\mathcal{E}^{(r)} \rightarrow \mathcal{E}^{(s)}$, for $r > s$; this map of principal bundles is $\rho^r_s$-equivariant.  We have a measurable trivialization $\mathcal{E}^{(r)} \cong M \times (\BR^n \times \GL^{(r)}(n))$.  

The fibers of $\pi^r_0$, each equivariantly diffeomorphic to $\GL^{(r)}(n)$, come equipped with a $\mathfrak{gl}^{(r)}(n)$-valued 1-form, which we will denote $\omega$.  For $g \in \GL^{(r)}(n)$, this form satisfies $R_g^* \omega = \Ad g^{-1} \circ \omega$, where $R_g$ is the right translation, and $\Ad$ is the adjoint representation. 

The \emph{$r$th prolongation of $F$} reflects the $\{ F^t_x \}$-action on $\mathcal{E}^{(r)}$ but is defined so that it leaves invariant the section $({\bf 0},\mbox{Id})$:
\begin{eqnarray*}
F^{(r)} & : & \BR \times M \rightarrow \mbox{Diff}^{q-r}(\BR^n \times \GL^{(r)}(n)) \\
 F^{(r)}(t,x) & : & (u,g) \mapsto (F^t_x(u), (J_u^{(r)} F^t_x) \cdot g \cdot (J_{\bf 0}^{(r)} F^t_x)^{-1})
\end{eqnarray*}

Here $J^{(r)}_u f$ denotes the $r$-jet of $f$ at $u$, identified with an element of $\GL^{(r)}(n)$ via pre- and post-composition with appropriate translations; thus the product in the second coordinate takes place in $\GL^{(r)}(n)$.
It is easy to verify that $F^{(r)}$ is a cocycle.
%$$ F^{(r)}(s+t,x)  = 
% (F^s_{\varphi^t(x)} \circ F^t_x(u), D_{F^t_x(u)}^{(k)} F^s_{\varphi^t(x)} D_u^{(k)} F^t_x g (D_u^{(k)} F^t_x)^{-1} (D_{F^t_x(u)}^{(k)} F^s_{\varphi^t(x)})^{-1} = 
%F^{(r)}(s,\varphi^t(x)) \circ F^{(r)}(t,x) $$ 

There is a corresponding \emph{$r$th prolongation of $T$}: 
$$T^{(r)}(t,x) = D_{({\bf 0},Id)} (F^{(r)}(t,x)) = D_{\bf 0} F^t_x \oplus \mbox{Ad}(J_{\bf 0}^{(r)} F^t_x)$$ 

In terms of $\omega$, the derivative of the $F^{(r)}(t,x)$-action on fibers of $\pi^r_0$ is by $\Ad J_{\bf 0}^{(r)} F^t_x$: If $v \in T_{(u,g)} \mathcal{E}^{(r)}_x \cap \ker (\pi^r_0)_*$ with $\omega(v) = X$, then we can write $v = ({\bf 0},gX)$.  Under the derivative of $F^{(r)}(t,x)$, the first coordinate remains ${\bf 0}$, and the second coordinate becomes
$$ (J_u^{(r)} F^t_x) gX (J_{\bf 0}^{(r)} F^t_x)^{-1} = (J_u^{(r)} F^t_x) g (J_{\bf 0}^{(r)} F^t_x)^{-1}  (J_{\bf 0}^{(r)}F^t_x) X (J_{\bf 0}^{(r)} F^t_x)^{-1}$$
which evaluates under $\omega$ to $(\Ad J_{\bf 0}^{(r)}F^t_x) (X)$.

\subsection{Subresonance Polynomials}
\label{sec.subres.poly}

Subresonance polynomials arise naturally in our context.  We will refer to \cite[Sec 3]{feres.normal.forms}; see also \cite[Prop 1.1]{guysinsky.katok.normal}.  An important role will be played by a nilpotent proper subgroup of the subresonance polynomials, to be defined below as the \emph{strict subresonance polynomials}.

Recall that $\lambda^{(s)} <0$ denotes the greatest element of the Lyapunov spectrum $\Sigma^0$ of the cocycle $\{ T^t_x \}$.  Let $d_i, i = 1, \ldots, s,$ be the dimensions of the Lyapunov distributions corresponding to $\lambda^{(1)}, \ldots, \lambda^{(s)}$, respectively; these dimensions are constant on $M_T$ by ergodicity of $\{\varphi^t \}$.  Let $W^1 \oplus \cdots \oplus W^s$ be a decomposition of $\BR^n$ into subspaces of dimension $d_1, \ldots, d_s$, respectively.  

We first establish some notation.  Given $l,r \geq 1$ and $\sigma \in \BR$, 
\begin{eqnarray*}
\widehat{W}^\sigma_l & := & \sum \{ (W^{p_1*} \otimes \cdots \otimes W^{p_l*} ) \otimes W^i \ : \lambda^{(i)} - \sum_{j=1}^l \lambda^{(p_j)} = \sigma \}  \\
\widehat{W}^{(r),\sigma} & := & \bigoplus_{l=1}^r \widehat{W}^{\sigma}_l; \qquad \widehat{V}^{\sigma}_l := \bigoplus_{\lambda \leq \sigma} \widehat{W}^\lambda_l; \qquad \widehat{V}^{(r),\sigma} = \bigoplus_{l=1}^r \widehat{V}^{\sigma}_l
\end{eqnarray*}

Let 
$$\Lambda = e^{\lambda^{(1)}} \id_{W_1} \oplus \cdots \oplus e^{\lambda^{(s)}} \id_{W_s} \in \GL(n).$$  

% \begin{definition}
% A \emph{subresonance} of $\Sigma^0$ of length $l$ is a multi-index $(p_1, \ldots, p_l, i)$ such that
% $$\lambda^{(i)} - \sum_{j=1}^{l} \lambda^{(p_j)} \leq 0,$$
% where $\lambda^{(i)}, \lambda^{(p_j)} \in \Sigma^0$.  
% A \emph{strict subresonance} of $\Sigma^0$ of length $l$ is a subresonance with
% $$\lambda^{(i)} - \sum_{j=1}^{l} \lambda^{(p_j)} \leq \lambda^{(s)}.$$
% A \emph{neutral subresonance} of $\Sigma^0$ of length $l$ is a subresonance with 
% $$\lambda^{(i)} - \sum_{j=1}^{l} \lambda^{(p_j)} = 0.$$

% The set of subresonances, strict subresonances, or neutral subresonances, respectively, of length $l$ will be denoted $\eta_l, \xi_l$, or $\eta_{l,0}$, respectively.  
% We will write
% $$ \eta^r = \bigcup_{l=1}^r \eta_l \qquad \xi^r = \bigcup_{l=1}^r \xi_l \qquad \eta^{r,0} = \bigcup_{l=1}^r \eta_{l,0}$$
% \end{definition}

\begin{definition}
\label{def.resonance.polys}
The \emph{subresonance polynomials} in $\GL^{(r)}(n)$ are 
$$ \mathcal{H}^{(r)}(\Lambda) =  \widehat{V}^{(r),0}  \cap \GL^{(r)}(n)$$
The corresponding \emph{strict subresonance polynomials} and \emph{resonance polynomials} are, respectively, 
$$  \mathcal{X}^{(r)}(\Lambda) = \id +  \widehat{V}^{(r),\lambda^{(s)}} \qquad \mbox{and} \qquad \mathcal{H}^{(r),0}(\Lambda) = \widehat{W}^{(r),0} \cap \GL^{(r)}(n)$$
%We set $\mathcal{X}^{(0)} = \{ \id \}$.
\end{definition}

Often $\mathcal{H}^{(r)}(\Lambda)$ will simply be written $\mathcal{H}^{(r)}$ when $\Sigma^0$ and $\Lambda$ are clear, and similarly for $\mathcal{X}^{(r)}$ and $\mathcal{H}^{(r),0}$.  
%Note that, with $\Lambda$ fixed, $\mathcal{H}^{(1)} \subset \cdots \subset \mathcal{H}^{(r)}$, and similarly for $\mathcal{X}^{(i)}$ and $\mathcal{H}^{(i),0}$.

Note that subresonance polynomials have bounded degree: for $r \geq \lfloor \lambda^{(1)}/\lambda^{(s)} \rfloor$, the groups $\mathcal{H}^{(r)} = \mathcal{H}^{(r+1)}$.  Note also that all derivatives at any $u \in \BR^n$ of a subresonance polynomial are strictly subresonance.  Similarly, the translation by $u \in \BR^n$ of a subresonance polynomial, $\tau_{- h(u)} \circ h \circ \tau_u$, is subresonance.

\begin{proposition}%(cf. \cite[Lemm 7, Prop 8]{feres.normal.forms})
\label{prop.subres.subgroups}
The polynomial germs in $\mathcal{H}^{(r)}, \mathcal{X}^{(r)}$, and $\mathcal{H}^{(r),0}$ form subgroups of $\GL^{(r)}(n)$, with Lie algebras
$$\mathfrak{h}^{(r)} = \widehat{V}^{(r),0} \qquad \mathfrak{x}^{(r)} = \widehat{V}^{(r),\lambda^{(s)}} \qquad \mathfrak{h}^{(r),0} = \widehat{W}^{(r),0}$$
The group $\mathcal{X}^{(r)}$ is a nilpotent normal subgroup of $\mathcal{H}^{(r)}$.  
When $r \geq \lfloor \lambda^{(1)}/\lambda^{(s)} \rfloor$, then polynomials in $\mathcal{H}^{(r)}$ are global diffeomorphisms of $\BR^n$.
% and $\mathcal{H}^{(r),0} \leq \mathcal{H}^{(r)}$.  
\end{proposition}

\begin{Pf}
It suffices to prove the statements for $r \geq \lfloor \lambda^{(1)}/\lambda^{(s)} \rfloor$, as the homomorphism $\rho^{r}_l$, $l \leq r$, preserves all group-theoretic properties and maps $\widehat{W}^{(r),\sigma}$ to $\widehat{W}^{(l),\sigma}$ and $\widehat{V}^{(r),\sigma}$ to $\widehat{V}^{(l),\sigma}$.

The \emph{$\Lambda^k$-stable polynomials} in $\alggl^{(r)}(n)$ are those $X$ for which $\interleave \Lambda^k X \Lambda^{-k} \interleave$ is bounded as $k \rightarrow \infty$.  The $\Lambda^k$-stable polynomials with invertible first derivative at ${\bf 0}$ form a subgroup of $\mbox{Diff}(\BR^n,{\bf 0})$ by \cite[Prop 8]{feres.normal.forms}; by \cite[Lemma 7 (6)]{feres.normal.forms}, this group equals $\mathcal{H}^{(r)}$.  The Lie algebra $\mathfrak{h}^{(r)}$ comprises all $\Lambda^k$-stable polynomials in $\alggl^{(r)}(n)$, which are $\widehat{V}^{(r),0}$.

We can express $\mathcal{X}^{(r)} \subset \mathcal{H}^{(r)}$ as
\begin{eqnarray*}
\mathcal{X}^{(r)} & = & \{ g \in \mathcal{H}^{(r)} \ : \ \lim_{k \rightarrow \infty} \frac{1}{k} \ln \interleave \Lambda^k g \Lambda^{-k} - \id \interleave \leq \lambda^{(s)} \} 
\end{eqnarray*}

By bilinearity of the bracket,
\begin{enumerate}
\item For $X,Y \in \sum_{\sigma \leq \lambda^{(s)}} \widehat{W}^{(r),\sigma} = \widehat{V}^{(r),\lambda^{(s)}}$,  
$$\lim_{k \rightarrow \infty} \frac{1}{k} \ln \interleave \Lambda^k [X,Y] \Lambda^{-k} \interleave  \leq \lim_{k \rightarrow \infty} \frac{1}{k} \left( \ln \interleave \Lambda^k X \Lambda^{-k} \interleave + \ln \interleave \Lambda^k Y \Lambda^{-k} \interleave \right)$$
\end{enumerate}

Using proposition \ref{prop.poly.bounds}, one can check the following inequalities:
\begin{enumerate}
\item[{(2)}] For $g, h \in \mathcal{H}^{(r)}$, 
$$ \lim_{k \rightarrow \infty} \frac{1}{k} \ln \interleave \Lambda^k g h \Lambda^{-k} - \id \interleave \leq  \lim_{k \rightarrow \infty} \frac{1}{k} \sup \{ \ln \interleave \Lambda^k g \Lambda^{-k} - \id \interleave, \ln \interleave \Lambda^k h \Lambda^{-k} - \id \interleave \}$$
% \item[($1'$)]  For $g, h \in \mathcal{H}^{(r)}$, 
% $$ \lim_{k \rightarrow \infty} \frac{1}{k} \ln \Vert \Lambda^k g h \Lambda^{-k} \Vert \leq  \lim_{k \rightarrow \infty} \frac{1}{k} \left( \ln \Vert \Lambda^k g \Lambda^{-k} \Vert + \ln \Vert \Lambda^k h \Lambda^{-k} \Vert \right)$$

\item[{(3)}] For $X \in \widehat{V}^{(r),\lambda^{(s)}}$ and $g \in \mathcal{H}^{(r)}$,
$$ \lim_{k \rightarrow \infty} \frac{1}{k} \ln \interleave \Lambda^k g X g^{-1} \Lambda^{-k} \interleave  \leq  \lim_{k \rightarrow \infty} \frac{1}{k} \ln \interleave \Lambda^k X \Lambda^{-k} \interleave$$
\end{enumerate}

Item (2) with $g,h \in \mathcal{X}^{(r)}$ shows this subset is closed under multiplication.  Similarly,
$$ \interleave \Lambda^k h^{-1} \Lambda^{-k} - \id \interleave = \interleave (\id - \Lambda^k h \Lambda^{-k}) (\Lambda^k h^{-1} \Lambda^{-k}) \interleave$$
and $\interleave \Lambda^k h^{-1} \Lambda^{-k} \interleave$ is bounded, so $\mathcal{X}^{(r)}$ is closed under inversion. Now $\mathcal{X}^{(r)}$ is a subgroup, and it is clear that the Lie algebra $\mathfrak{x}^{(r)}$ is as claimed.  

One deduces from item (1) that $\mathfrak{x}^{(r)}$ is nilpotent.  

Finally, from (3), $\mathfrak{x}^{(r)}$ is a Lie algebra ideal, so $\mathcal{X}^{(r)} \lhd \mathcal{H}^{(r)}$.

We will prove inductively on $r$ that $\mathcal{H}^{(r),0}$ forms a subgroup (dropping the assumption $r \geq \lfloor \lambda^{(s)}/\lambda^{(1)} \rfloor$).  It is easy to see that (truncations of) compositions of polynomials in $\widehat{W}^{(r),0}$ are again in $\widehat{W}^{(r),0}$, so we will just check closure under inverse.  

For $r=1$, we have
$$\mathcal{H}^{(1),0} = \GL(W^1) \oplus \cdots \oplus \GL(W^s)$$
with Lie algebra
$$\mathfrak{h}^{(1),0} = \End(W^1) \oplus \cdots \oplus \End(W^s) = \widehat{W}^{(1),0}$$

Assume that $\mathcal{H}^{(j-1),0}$ is closed under inversion.  Let $g,h \in \mathcal{H}^{(j)}$ with $\overline{g} = \rho^j_{j-1}(g), \overline{h} = \rho^j_{j-1}(h) \in \mathcal{H}^{(j-1),0}$.  Write $g = \overline{g} + \gamma$ and $h = \overline{h} + \delta$ with $\gamma, \delta \in S^{(j)}(n)$.  Then the composition in $\mathcal{H}^{(j)}$ can be written
\begin{eqnarray}
\label{eqn.Hr0.composition}
g \circ h & = & \overline{g} \circ \overline{h} + (D_{\bf 0} g) \circ \delta + \gamma \circ (D_{\bf 0} h) + R
\end{eqnarray}
where  $\overline{g} \circ \overline{h} \in \mathcal{H}^{(j-1),0}$, and $R$ is a sum of degree-$j$ compositions of terms of $\overline{g}$ and $\overline{h}$. Assume $g \in \mathcal{H}^{(j),0}$, and let $h = g^{-1} = \overline{g}^{-1} + \delta$.  Then
\begin{eqnarray}
\label{eqn.Hr0.inverse}
\delta = - D_{\bf 0} g^{-1} \circ \gamma \circ D_{\bf 0} g^{-1} - D_{\bf 0} g^{-1} \circ R \in \widehat{W}^0_j   
\end{eqnarray}
We conclude that $\mathcal{H}^{(r),0} \leq \mathcal{H}^{(r)}$ with Lie algebra $\mathfrak{h}^{(r),0} = \widehat{W}^{(r),0}$.  
\end{Pf}

\begin{proposition}(compare \cite[Prop 7]{feres.normal.forms})
\label{prop.subres.algebraic}
Let $\{F_k \}$ be a sequence in $\GL^{(r)}(n)$, $r \geq 2$, with $\rho^r_1(F_k) = T_k$.  Let $F^k = F_k \circ \cdots \circ F_1$ and $T^k = T_k \circ \cdots \circ T_1$.  Suppose that $(T^{k*} T^k)^{1/2k} \rightarrow \Lambda$, with spectrum $\Sigma^0$ and eigenspace decomposition $W^1 \oplus \cdots \oplus W^s$.  Then:
\begin{enumerate}
% \item $$\lim \sup_m \frac{1}{m} \ln \Vert F^m \cdot X \cdot (F^m)^{-1} \Vert \leq \lambda^{(s)}$$
\item Denote by $A^k$ the restriction of $\mbox{Ad} (F^k)$ to $S^{(r)}(n)$.  The limit $\lim_{k \rightarrow \infty}$ $(A^{k*}A^k)^{1/2k}$ exists, and has spectrum comprising all 
$$\exp(\lambda^{(i)} - \sum_{j=1}^r \lambda^{(p_j)}) \ \qquad \mbox{where} \qquad  \ \lambda^{(i)}, \lambda^{(p_j)} \in \Sigma^0$$
\item For $X \in S^{(r)}(n)$, the limit $\lim_{k \rightarrow \infty} \frac{1}{k} \ln \interleave F^k \circ X \circ F^{-k} \interleave$
exists, and equals $\sigma$ if and only if $X \in \widehat{W}_{r}^\sigma$.
\end{enumerate}
\end{proposition}

\begin{Pf}
For $X \in S^{(r)}(n)$ and $F \in \GL^{(r)}(n)$ with $\rho^r_1(F) = T$, the conjugation is
$$ F \circ X \circ F^{-1} =  T \circ X \circ T^{-1}  \qquad \mbox{(see \cite[Eqn. 27]{feres.normal.forms})}$$ 

Then 
$$ (A^{*k} A^k)^{1/2k} (X) = ((\Ad T^k)^* (\Ad T^k))^{1/2k}(X)$$

The symmetrization of $\Ad T^k$ on $S^{(r)}(n)$ is simply $\mbox{Ad}(T^{k*} T^k)^{1/2}$.  Then the above expression is  
$$(T^{k*} T^k)^{1/2k} \circ X \circ (T^{k*} T^{k})^{-1/2k} \rightarrow \Lambda \circ X \circ \Lambda^{-1}$$
The remaining claims of (1) and (2) now follow.
\end{Pf}

\subsection{$T^{(r)}$ spectrum, Lyapunov decomposition, and algebraic hull}
\label{subsec.Tr.spectrum}

Let $T^{(r)}(t,x) = D_{({\bf 0}, \id)} F^{(r)}(t,x)$ as in section \ref{sec.prolongation} above.  We have assumed that $q \geq 2$ and 
$$\sup_{0 \leq \epsilon \leq 1} \ln^+ \Vert F^{\pm \epsilon}_x \Vert_q \in L^1(M,\mu).$$  
We will deduce that, for $r \leq q$, the linear cocycle $T^{(r)}$ satisfies (\ref{eqn.met}). 
%$$\sup_{0 \leq \epsilon \leq 1} \ln^+ \Vert T^{(r)}(\epsilon,x) \Vert \in L^1(M,\mu).$$ 
The two components of $T^{(r)}(t,x)$ are $D_{\bf 0} F^t_x$ and $\Ad J_{\bf 0}^{(r)} F^t_x$.  As $\Vert D_{\bf 0} F^t_x \Vert \leq \Vert F^t_x \Vert_q$, we can focus on the second factor.  
%We will use the norm defined on $r$-jets in section \ref{sec.prolongation}.  
Fix $x \in M$ and $\epsilon$ between $-1$ and $1$, and let $X \in \mathfrak{gl}^{(r)}(n)$ with $\interleave X \interleave = 1$.  By proposition \ref{prop.poly.bounds},
\begin{eqnarray*}
\ln^+ \interleave (\Ad J_{\bf 0}^{(r)} F^\epsilon_x) (X) \interleave & \leq & \ln^+ \left( \interleave J_{\bf 0}^{(r)} F^\epsilon_x \interleave \cdot \max_{i=1,r} \interleave X \circ (J_{\bf 0}^{(r)} F^\epsilon_x)^{-1} \interleave^i \cdot c(r,n) \right) \\
% & \leq & \ln^+ \Vert J_{\bf 0}^{(r)} F^\epsilon_x \Vert + \ln^+ \left( C(r,n) \sup_{1 \leq i \leq r} \Vert X \circ (J_{\bf 0}^{(r)} F^\epsilon_x)^{-1} \Vert^i \right) \\
& \leq &  \ln^+ \interleave J_{\bf 0}^{(r)} F^\epsilon_x \interleave + \ln^+ \interleave (J_{\bf 0}^{(r)} F^\epsilon_x)^{-1} \interleave^{r^2}  + \ln c(r,n)^{r+1} 
%& \leq &  \ln^+ \Vert J_{\bf 0}^{(r)} F^\epsilon_x \Vert + \ln^+  \Vert (J_{\bf 0}^{(r)} F^\epsilon_x)^{-1} \Vert^r + C_0
\end{eqnarray*}

Let $C_0 = (r+1) \ln c(r,n)$.  As the above expression is independent of the unit vector $X$, we obtain with respect to the linear norm on $\GL(\alggl^{(r)}(n))$,  
\begin{eqnarray*}
& & \int \ln^+ || (\Ad J_{\bf 0}^{(r)} F^\epsilon_x) || dx  \\
& \leq & \int \ln^+ \interleave  J_{\bf 0}^{(r)} F^\epsilon_x \interleave dx + \int r^2 \ln^+ \interleave (J^{(r)}_{\bf 0} F^\epsilon_x)^{-1} \interleave dx + C_0 \\
& \leq &  \int \ln^+ \interleave  J_{\bf 0}^{(r)} F^\epsilon_x \interleave dx + \int r^2 \ln^+ \interleave J^{(r)}_{\bf 0} F^{-\epsilon}_x \interleave dx + C_0 \\
& \leq &  \int \ln^+ \Vert F^\epsilon_x \Vert_r dx + \int r^2 \ln^+ \Vert F^{-\epsilon}_x \Vert_r dx + C_0
\end{eqnarray*}

%These two integrals exist by our original assumption, because $r \leq q$. It is now straightforward to conclude that,
Inserting $\sup_\epsilon$ in the integrands in the above chain of inequalities leads to
$$\sup_{0 \leq \epsilon \leq 1} \ln^+ \Vert \Ad J_{\bf 0}^{(r)} F^{\pm \epsilon}_x \Vert \in L^1(M,\mu),$$
as desired.  The Multiplicative Ergodic Theorem thus applies to $\{ T^{(r)}(t,x) \}$.

Denote by $M_T \subseteq M$ the set on which the conclusions of the MET hold for $\{T^{(r)}(t,x)\}$.  For $\sigma \in \mbox{Spec } T^{(r)}$, denote by $V_x^{(r),\sigma}$ and $W^{(r),\sigma}_x$ the corresponding terms in the Lyapunov filtration and Lyapunov decomposition, respectively.  We can assume, after a tempered linear cocycle equivalence, that the Lyapunov decomposition for $\{ T^t_x \}$ is constant $W_x^i = W^i$, where $W^i$ is as in section \ref{sec.subres.poly} (see \cite[Thm 6.1]{barreira.pesin.introset}).  The following proposition implies that now the Lyapunov decomposition for $\{ T^{(r)}(t,x) \}$ is also constant.

%% Denote by $\underline{\sigma}$ the maximal element of $\mbox{Spec } T^{(r)}$ less than $\sigma$, and $\overline{\sigma}$ the minimal element greater than $\sigma$, so $V_x^{(r),\underline{\sigma}} \subset V_x^{(r),\sigma} \subset V_x^{(r),\overline{\sigma}}$.

\begin{proposition}
\label{prop.Tr.spectrum}
Let $1 \leq r \leq q$.  
\begin{enumerate}
\item The Lyapunov exponents of $\{ T^{(r)}(t,x) \}$ are $\cup_{l=0}^r \Sigma^l$, where,  
$$ \Sigma^l = \{ \lambda^{(i)} - \sum_{j=1}^l \lambda^{(p_j)} \ : \ \lambda^{(i)}, \lambda^{(p_j)} \in \Sigma^0 \} \qquad \mbox{for} \ l \geq 1$$
\item  For $\sigma \in \cup_{l=0}^r \Sigma^l$, for $x$ in a $\varphi^t$-invariant $\mu$-conull set, the $\{ T^{(r)}(t,x) \}$ Lyapunov space at $x$ comprises $W^\sigma \subset \BR^n$, direct sum
$$ \widehat{W}^{(r),\sigma} = W_x^{(r),\sigma} \cap \alggl^{(r)}(n)$$ 
\end{enumerate}  
\end{proposition}

\begin{Pf}
For $r \geq 2$, proposition \ref{prop.subres.algebraic} says that $\Sigma^{r} \subset \mbox{Spec } T^{(r)}$ and that, for each $\sigma \in \Sigma^r$,  
$$ W^{(r),\sigma}_x \cap S^{(r)}(n) = \widehat{W}^\sigma_r.$$  

Thus for almost all $x \in M$, by proposition \ref{prop.basic.lyapunov} (1) and proposition \ref{prop.subres.algebraic} (1),
$$ S^{(r)}(n) = \bigoplus_{\sigma \in \Sigma^r} \widehat{W}^\sigma_r$$

The quotient vector bundle is 
$$ T_{({\bf 0}, \id)} \mathcal{E}^{(r-1)} = T_{({\bf 0},\id)} \mathcal{E}^{(r)} / S^{(r)}(n)$$
The spectrum of $T^{(r-1)}$ is a subset of $\mbox{Spec } T^{(r)}$, and $W^{(r),\sigma}_x$ projects onto $W^{(r-1),\sigma}_x$, by proposition \ref{prop.basic.lyapunov} (2).  Projecting $\mathcal{E}^{(j)} \rightarrow \mathcal{E}^{(j-1)}$ for $j=r, \ldots, 2$, yields $\mbox{Spec } T^{(r)} = \cup_{l=2}^r \Sigma^l \cup \mbox{Spec } T^{(1)}$.  It is easy to see that $\mbox{Spec } T^{(1)} = \Sigma^1 \cup \Sigma^0$, as the action here only involves $T^t_x$ and $\Ad T^t_x$; moreover, the Lyapunov space 
$$W^{(1),\sigma}_x = W^\sigma \oplus \widehat{W}^\sigma_1.$$
We conclude $\mbox{Spec } T^{(r)} = \cup_{l=0}^r \Sigma^l$, as claimed in (1).

With respect to the filtration induced by $\ker (\rho^r_j)_*$, $j=r-1,\ldots, 0,$ the space $W^{(r),\sigma}_x \cap \alggl^{(r)}(n)$ has the same associated graded space  as $\widehat{W}^{(r),\sigma}$; in particular, $W^{(1),\sigma}_x \cap \alggl(n) = \widehat{W}^\sigma_1$.  Suppose that $W^{(j-1),\sigma}_x \cap \alggl^{(j-1)}(n) = \widehat{W}^{(j-1),\sigma}$ for some $j \geq 2$.  Form the tensor bundle
\begin{eqnarray*}
\mathcal{T}_x & = & (W^{(j-1),\sigma}_x \cap \ker(\rho^{j-1}_0)_*)^* \otimes \left( \ker(\rho^{j}_{j-1})_* / (W^{(j),\sigma}_x \cap \ker(\rho^j_{j-1})_*) \right) \\
& = & (\widehat{W}^{(j-1),\sigma})^* \otimes (S^{(j)}(n) / \widehat{W}^\sigma_j)
\end{eqnarray*}

Let $\{T^{(j)}(t,x)_*\}$ be the linear cocycle on $\mathcal{T}$ induced from $\{ T^{(j)}(t,x)\}$, and define an invariant tensor $\{ \tau_x \in \mathcal{T}_x \}$ by
$$ W^{(j),\sigma}_x = \mbox{graph}(\tau_x) + \widehat{W}^\sigma_j$$
Proposition \ref{prop.T.stability} says that $\tau_x$ belongs, for almost every $x$, to the $\{T^{(j)}(t,x)_* \}$ Lyapunov space $\mathcal{T}^0_x$.  But this space is trivial, so $\tau_x = 0$, and $W^{(j),\sigma}_x \cap \alggl^{(j)}(n)= \widehat{W}^{(j),\sigma}$ for almost every $x$.  Induction on $j$ yields (2).   
\end{Pf}

Henceforth denote $W^{(r),\sigma} = W^\sigma \oplus \widehat{W}^{(r),\sigma}$ and similarly for $V^{(r),\sigma}$.

\begin{proposition}
\label{prop.J0.subres}
For all $x \in M_T$ and all $t \in \BR$,  
$$J^{(r)}_{\bf 0} F^t_x \in \mathcal{H}^{(r),0}.$$
If $(\{ G^k_x \},\psi) \in Z(F^t_x)$ then also for all $k$ and almost all $x$,
$$ J_{\bf 0}^{(r)} G^k_x \in \mathcal{H}^{(r),0}.$$
%\in \mbox{Diff}^q(\BR^n,{\bf 0})$ is a cocycle on $\BZ \times M$ satisfying ?*? for which the prolonged linear cocycle 
\end{proposition}

\begin{remark}
In fact, the conclusion holds for any $\{ G^k_x \}$ satisfying (1) and (3) of definition \ref{def.centralizer}, for which the prolonged linear cocycle
$$ U^{(r)}(k,x) = D_{\bf 0} G^k_x \oplus \mbox{Ad}(J^{(r)}_{\bf 0} G^k_x)$$
preserves the Lyapunov decomposition of $T^{(r)}$.
\end{remark}

\begin{remark}
\label{rmk.rjets.contract}
Viewing $J_{\bf 0}^{(r)} F^k_x$ as a polynomial map in $\GL^{(r)}(n)$, we see from remark \ref{rmk.jets.arb.small} that it tends to ${\bf 0}$ as $k \rightarrow \infty$. 
\end{remark}

\begin{Pf}
The proof only makes use of the Lyapunov decomposition of $\{ T^{(r)}(t,x) \}$; recall from proposition \ref{prop.Tr.spectrum} (2), 
$$W^{(r),0}_x \cap \alggl^{(r)}(n)= \mathfrak{h}^{(r),0}$$ 
The proof will apply essentially verbatim to $\{ G^k_x \}$.

First, $T^t_x = D_{\bf 0} F^t_x$ preserves the Lyapunov decomposition $W^1 \oplus \cdots \oplus W^s$, which implies that, for all $x \in M_T$,
\begin{eqnarray}
\label{eqn.J1.H10} 
J_{\bf 0}^{(1)} F^t_x \in \mathcal{H}^{(1),0} =  \GL(W^1) \times \cdots \times \GL(W^s)
\end{eqnarray}

For $j \geq 2$, let $g = J_{\bf 0}^{(j)} F^t_x$, and $\overline{g} = \rho^j_{j-1} (g)$.  Write $g = \overline{g} + \gamma$, with $\gamma \in S^{(j)}(n)$, and suppose $\overline{g} \in \mathcal{H}^{(j-1),0}$.
%  Now it easily follows that, if , then
% \begin{itemize}
% \item $ g \circ h \equiv D_{\bf 0} g \circ \delta + \gamma \circ D_{\bf 0} h$ mod $\mathfrak{h}^{(j),0}$. 
% \item $g^{-1} \equiv - D_{\bf 0} g^{-1} \circ \gamma \circ D_{\bf 0} g^{-1}$ mod $\mathfrak{h}^{(j),0}$;
% \end{itemize}
Let $\Lambda$ be as in section \ref{sec.subres.poly}, and note that $\Lambda \in \mathcal{H}^{(1),0} \subset \mathfrak{h}^{(j),0}$.  By (\ref{eqn.J1.H10}), $\Lambda$ commutes with $D_{\bf 0} F^t_x$ for all $t \in \BR$, $x \in M_T$.  By similar computations to equations (\ref{eqn.Hr0.composition}) and (\ref{eqn.Hr0.inverse}),
$$ g \circ \Lambda \circ g^{-1} \equiv \gamma \circ \Lambda \circ  D_{\bf 0} g^{-1} - \Lambda \circ \gamma \circ D_{\bf 0} g^{-1} \ \mbox{mod } \mathfrak{h}^{(j),0}$$

On the other hand, this conjugate belongs to $\mathfrak{h}^{(j),0}$, and both $\Lambda$ and $D_{\bf 0}g^{-1}$ preserve all the Lyapunov spaces $W^i$, so
$$ \Lambda \circ \gamma \circ \Lambda^{-1} \equiv \gamma \ \mbox{mod } \mathfrak{h}^{0}_j,$$
which implies $\gamma \in \mathfrak{h}^0_j$.  Proceeding inductively gives $J_{\bf 0}^{(r)} F^t_x \in \mathcal{H}^{(r),0}$ for all $t$ and all $x \in M_T$.  
\end{Pf}

\subsection{Dynamical submanifolds as reductions of $\mathcal{E}^{(r)}$}
\label{subsec.xr.reduction}

In this section we apply theorem \ref{thm.ruelle.submanifolds} to $\{ F^{(r)}(t,x) \}$ and interpret the resulting submanifolds as invariant $\mathcal{X}^{(r)}$-reductions of $\mathcal{E}^{(r)}$.  Combining this information with proposition \ref{prop.J0.subres} gives an $\mathcal{H}^{(r)}$-reduction of $\mathcal{E}^{(r)}$ invariant by the natural action of $\{ F^t_x \}$.

First we must recall an aspect of the construction of the $\intnu_x^\lambda$ that will play a role in our upcoming proof.

\subsubsection{Digression on Ruelle's proof: expressing as graphs} 
\label{subsubsec.digression.graph}
For $\lambda^{(i)} < \lambda < \lambda^{(i+1)}$, the submanifold $\intnu_x^\lambda$ is defined on \cite[p 47]{ruelle.diff.ergodic} as the image under a certain map $\Phi$ of the graph of a $C^1$ function 
$$ \varphi: V_x^i \cap B(\alpha(x)) \rightarrow (V_x^i)^\perp \cap B(\alpha(x))$$

Write $\alpha = \alpha(x)$.  The equation \cite[(5.15)]{ruelle.diff.ergodic} gives the bound
$$ \Vert D_u \varphi \Vert  \leq \frac{A \delta \sqrt{1-(A \delta)^2}}{1-2(A \delta)^2} \qquad \forall \ u \in V_x^i \cap B(\alpha)$$
where $A$ and $\delta$ are the constants given by the perturbation theorem as described in subsection \ref{subsubsec.digression}.  In equation (5.9) of \cite{ruelle.diff.ergodic}, the additional bound $\delta < 1/\sqrt{2} A$ was imposed.  There is no problem in assuming $\delta < 1/2A$, which makes $\Vert D_u \varphi \Vert < 1$.

The map $\Phi : (V_x^i \cap B(\alpha)) \oplus ((V_x^i)^\perp \cap B(\alpha)) \rightarrow B(\alpha)$ is
$$ \Phi(u, v) = \frac{u}{\alpha} \sqrt{\alpha^2 - \Vert v \Vert^2} + v$$
For a differentiable path $\gamma(t) = (u(t),v(t))$ in $\hat{L}_b = \Phi^{-1}(\{ u = b \})$, 
$$ \dot{u} = \frac{\langle v, \dot{v} \rangle}{\alpha^2 - \Vert v \Vert^2} \cdot u$$
If $\gamma$ is contained in $B(\alpha /\sqrt{2})$, then $\Vert \dot{v}(t) \Vert \geq \Vert \dot{u}(t) \Vert$ along $\gamma$ and $\Vert v(1) - v(0) \Vert \geq \Vert u(1) - u(0) \Vert$.  Now $\mbox{graph }\varphi$ intersects $\hat{L}_b \cap B(\alpha/\sqrt{2})$ in at most one point, and 
$$ \intnu_x^{\lambda} \cap B(\alpha/2) \subseteq \Phi(\mbox{graph }\varphi \cap B(\alpha/\sqrt{2})) \subseteq \intnu_x^\lambda \cap B(\alpha/\sqrt{2})$$ 
is also the graph of a $C^1$ function. End of digression.

\bigskip

A fiber $\mathcal{E}_x^{(r)} = \BR^n \times \GL^{(r)}(n)$ can be viewed as an open subset of $\BR^N$ with origin at $({\bf 0},\id)$.  Pulling back the norm $\interleave \cdot \interleave$ on $\alggl^{(r)}(n)$ by $\omega$ yields a Finsler metric on $\GL^{(r)}(n)$, and thus a product Finsler metric on each $\mathcal{E}^{(r)}_x$. This metric is comparable to the Euclidean metric on $\BR^N$ on a compact neighborhood of the origin. 
%, for which $\pi^r_0$ is a Riemannian submersion.

The verification that, for $r \leq q$, 
$$ \sup_{0 \leq \epsilon \leq 1} \ln^+ \Vert F^{(r)}(\pm \epsilon,x) \Vert_{q-r} \in L^1(M,\mu)$$
is left to the reader, with the indication that it resembles the verification at the beginning of section \ref{subsec.Tr.spectrum} that $\{ T^{(r)}(t,x) \}$ satisfies (\ref{eqn.met}).  Let $M_F$ be the set where the conclusion of the Ergodic Theorem holds for this function under $\{ \varphi^t \}$.  Let $\sigma < 0$ in $\Sigma^l$, $l \leq r$.  Theorem \ref{thm.ruelle.submanifolds} gives $C^{q-r}$-smooth submanifolds $\nu_x^{(r),\sigma} \subset \mathcal{E}^{(r)}_x$ for all $x \in M_0 = M_F \cap M_T$.  Denote by $\nu_x^{(r),s}$ the submanifolds associated to $\sigma = \lambda^{(s)}$.  The following key result is a nonlinear analogue of proposition \ref{prop.Tr.spectrum} (2).

\begin{proposition}
\label{prop.xr.reduction}
For all $x$ in a $\varphi^t$-invariant, $\mu$-conull subset of $M_0$ and $1 \leq r \leq q$, $\nu_x^{(r),s}$ is a $C^{q-r}$-smooth reduction of $\mathcal{F}^{(r)} \mathcal{E}_x$ to $\mathcal{X}^{(r)}$, projecting onto $\nu_x^s$ under $\pi^r_0$.
\end{proposition}

\begin{Pf}

\emph{Step 1: Fibers under $\pi^r_{r-1}$ tangent to $\mathfrak{x}^{(r)}$}.   
For $(w,v) \in T_{(u,g)} \mathcal{E}^{(r)}_x$
% \cong T_u\BR^n \oplus T_g (\GL^{(r)}(n)),$ 
the effect of the action of $F^{(r)}(t,x)$ in the Finsler metric $\interleave \cdot \interleave_{(u,g)}$ is
$$ \interleave D_{(u,g)}F^{(r)}(t,x)\cdot (w,v) \interleave^2 = \Vert D_u F^t_x(w) \Vert^2 + \interleave (\Ad J^{(r)}_{\bf 0} F^t_x)(\omega(v)) \interleave^2 $$
%Let $u \in \mathcal{E}_x$.  Use perturbation theorem to show that spectrum of $T'(k,x) = DF(k,x)(u)$ is the same as that of $T(k,x)$.  Then can apply p
If $(u,g) \in \nu_x^{(r),s}$, then $F^{(r)}(k,x) \cdot (u,g)$ is in the Euclidean ball $\overline{B(1)}$ for all sufficiently large $k$.  The characterization of $\nu^{(r),s}_x$ in theorem \ref{thm.ruelle.submanifolds} can equivalently be expressed in terms of our Finsler distance $d$. 

Let $r \geq 2$ and $v \in T_{(u,g)} \nu^{(r),s}_x \cap \ker (\pi^r_{r-1})_*$.  Then 
$$ \lim_{k \rightarrow \infty} \frac{1}{k} \ln \interleave F^{(r)}(k,x)_* v \interleave = \lim_{k \rightarrow \infty} \frac{1}{k} \ln \interleave ( \Ad J^{(r)}_{\bf 0} F^k_x ) \cdot (\omega(v)) \interleave \leq \lambda^{(s)}$$
which implies, by proposition \ref{prop.subres.algebraic} (2), 
%applied to $F^k = J^{(r)}_{\bf 0} F^k_x$, gives
$$ \omega(v) \in \widehat{V}^{\lambda^{(s)}}_r = \mathfrak{x}^{(r)} \cap S^{(r)}(n)$$
We denote this last space by $\widehat{V}^s_r$ below.
%note that it equals $\mathfrak{x}^{(r)} \cap S^{(r)}(n)$ by proposition \ref{prop.subres.subgroups}.  
The above implication is also clear when $r=1$, with the conclusion $\omega(v) \in \widehat{V}_1^s$.

Conversely, if $X \in \widehat{V}^s_r$ with $\interleave X \interleave = 1$, then, by proposition \ref{prop.subres.algebraic} (2), 
$$ \lim_{k \rightarrow \infty} \frac{1}{k} \ln \interleave (J_{\bf 0}^{(r)} F^k_x) \circ X \circ (J_{\bf 0}^{(r)} F^{k}_x)^{-1} \interleave \leq \lambda^{(s)}.$$
Again, the same bound holds in the case $r=1$.
For $g' = g \cdot e^{\eta X}$, $\eta \in \BR$,
\begin{eqnarray*}
d(F^{(r)}(k,x) \cdot (u,g'), F^{(r)}(k,x) \cdot (u,g) ) & \leq & |\eta| \cdot \interleave D_{(u,h)} F^{(r)}(k,x) \cdot (0,\omega^{-1}(X)) \interleave \\
& = & | \eta| \cdot \interleave (\Ad J_{\bf 0}^{(r)} F^k_x)(X) \interleave 
\end{eqnarray*}
where $h$ is some point in $\{ g e^{t \eta X} \}_{t =0}^1$. 
%% Then
%% $$ \lim_{k \rightarrow \infty} \frac{1}{k} \ln d\left( F^{(r)}(k,x) \cdot (u,g'), F^{(r)}(k,x) \cdot (u,g) \right) \leq \lambda^{(s)}$$
As $(u,g)$ belongs to $\nu^{(r),s}_x$,  
 $$ \lim_{k \rightarrow \infty} \frac{1}{k} \ln d\left( F^{(r)}(k,x) \cdot (u,g'), ({\bf 0}, \id) \right) \leq \lambda^{(s)}$$ 
%(Here the norm denotes the distance to $({\bf 0}, \id)$.)  
and $(u,g') \in \nu^{(r),s}_x$ for  $\eta \in \BR$ (see proposition \ref{prop.submanifolds.sequence}).  Thus for all $r$,
$$\widehat{V}^s_r =  \omega(T_{(u,g)} \nu^{(r),s}_x \cap \ker(\pi^r_{r-1})_*).$$
Moreover, $\nu^{(r),s}_x$ is saturated by the principal action of $\mathcal{X}^{(r)} \cap \ker \rho^r_{r-1}$.  

\medskip

\emph{Step 2: Projection under $\pi^r_{r-1}$ is onto}. 

The image $\pi^r_{r-1}(\nu_x^{(r),s})$ is clearly contained in $\nu_x^{(r-1),s}$.  
%We have a $T^{(r)}$-invariant subbundle
%$$ \widehat{V}^s_r \subseteq T_{({\bf 0}, \id)} \nu_x^{(r),s} = V^{(r),s}_x$$ 
%\subset T_{({\bf 0},\id)} \mathcal{E}^{(r)}$$ 
%= \oplus_{\sigma \leq \lambda^{(s)}} W_x^{(r),\sigma}$$
The image of $T_{({\bf 0}, \id)} \nu_x^{(r),s}$ under $(\pi^r_{r-1})_{*}$ is, by proposition \ref{prop.Tr.spectrum} (2) and theorem \ref{thm.ruelle.submanifolds},
$$V^{(r),s}/\widehat{V}^s_r = V^{(r-1),s} = T_{({\bf 0},\id)}\nu_x^{(r-1),s}$$
% By proposition \ref{prop.Tr.spectrum} (2), this space projects under $(\pi^r_{r-1})_*$ onto $T_{({\bf 0}, \id)} \nu_x^{(r-1),s} = V^{(r-1),s}_x$.  
Thus for $x$ in a $\varphi^t$-invariant, conull set, $\pi^r_{r-1}(\nu^{(r),s}_x)$ and $\nu_x^{(r-1),s}$ coincide in a neighborhood of $({\bf 0},\id)$.   

Choose $\lambda$ between $\lambda^{(s)}$ and the next greater element of $\mbox{Spec } T^{(r)}$; choose $\lambda/3 < \zeta < 0$.  We claim that for all $k \geq 0$, there exists $m \geq 0$ such that
$$ F^{(r-1)}(m,\varphi^k(x)) \cdot \intnu_{\varphi^k(x)}^{(r-1),\lambda}  \subset \pi^r_{r-1} \cdot \nu_{\varphi^{k+m}(x)}^{(r),s}$$
The intersection of these contains a neighborhood in each of $({\bf 0}, \id)$.  
From the proof of proposition \ref{prop.submanifolds.sequence}, the radii of the left-hand terms are bounded above by $\beta(\varphi^k(x)) e^{m \lambda}$.  Also, $\nu_{\varphi^{k+m}(x)}^{(r),s}$ contains a ball of radius $\alpha(\varphi^{k+m}(x)) > C e^{(k+m)\zeta}$; as the restriction of $\pi^r_{r-1}$ is the quotient by the proper action of $\mathcal{X}^{(r)} \cap \ker \rho^r_{r-1}$, the projections %$\pi^r_{r-1} \cdot \nu_{\varphi^{k+m}(x)}^{(r),s}$ 
contain balls of comparable radius.  The claim follows.  Finally,
\begin{eqnarray*}
\pi^r_{r-1} \cdot \nu_x^{(r),s} & = & \pi^r_{r-1} \cdot \bigcup_{k \geq 0} F^{(r)}(k,x)^{-1} \cdot \intnu_{\varphi^k(x)}^{(r),\lambda} \\
& = & \pi^r_{r-1} \cdot \bigcup_{k \geq 0} F^{(r)}(k,x)^{-1} \cdot \nu_{\varphi^k(x)}^{(r),s} \\
& = & \bigcup_{k \geq 0} F^{(r-1)}(k,x)^{-1} \cdot \pi^r_{r-1} \cdot \nu_{\varphi^k(x)}^{(r),s} \\
& \supset & \bigcup_{k \geq 0}  F^{(r-1)}(k,x)^{-1} \cdot \intnu_{\varphi^k(x)}^{(r-1),\lambda} \\
& = & \nu_x^{(r-1),s}
\end{eqnarray*}

\medskip

\emph{Step 3: Fibers under $\pi^r_0$ tangent to $\mathfrak{x}^{(r)}$}.

Now $\pi^r_0 = \pi^1_0 \circ \cdots \circ \pi^r_{r-1}$ maps $\nu^{(r),s}_x$ onto $\nu^s_x$.  The vertical tangent subspaces of $\nu^{(r),s}_x$ evaluate under $\omega$ to subspaces 
$$\mathfrak{y}_{x}^{(r)}(u,g) = \omega(T_{(u,g)} \nu^{(r),s}_x \cap \ker (\pi^r_0)_*)\subset \alggl^{(r)}(n)$$
with the same associated graded algebra as $\mathfrak{x}^{(r)}$.  
%To finish the proof, it suffices to show that these subspaces everywhere equal $\mathfrak{x}^{(r)}$.
Note that
$$ \mathfrak{y}_{\varphi^t(x)}^{(j)} (F^{(j)}(t,x) \cdot (u,g) ) = (\Ad J^{(j)}_{\bf 0} F^t_x) \cdot \mathfrak{y}_x^{(j)}(u,g)$$
for $1 \leq j \leq r$, and they vary smoothly in $(u,g)$ and measurably in $x$.

Suppose that $\mathfrak{y}_x^{(j-1)}(u,g) = \mathfrak{x}^{(j-1)}$ for almost all $x$, for all $u,g$.  
%Note that $\mathfrak{y}_x^{(j-1)}(u,g) = (\rho^j_{j-1})_*(\mathfrak{y}_x^{(j)}(u,g))$.  
As $\mathfrak{y}_x^{(j)}(u,g) \cap S^{(j)}(n) = \widehat{V}^s_j$ from step 1, we can define linear functions 
$$\tau_x(u,g) : \mathfrak{x}^{(j-1)} \rightarrow S^j(n)/\widehat{V}^s_j \ \ \mbox{by} \ \ \mathfrak{y}_x^{(j)}(u,g)= \mbox{graph}(\tau_x(u,g)) + \widehat{V}^s_j$$
%% Also $\tau$ varies smoothly in $(u,g)$ and measurably in $x$.  We have
%% $$ \mathfrak{y}_{\varphi^t(x)}^{(j)} (F^{(j)}(t,x) \cdot (u,g) ) = (\Ad J^{(j)}_{\bf 0} F^t_x) \cdot \mathfrak{y}_x^{(j)}(u,g)$$
Note that $\tau$ corresponds to a tensor on $\cup_x \nu_x^{(r),s}$, smooth in $(u,g)$ and $F^{(j)}$-invariant---that is,
$$ \tau_{\varphi^t(x)}(F^{(j)}(t,x)\cdot (u,g)) = (\Ad J^{(j)}_{\bf 0} F^t_x) \circ \tau_x(u,g) \circ (\Ad J^{(j-1)}_{\bf 0} F^t_x)^{-1}$$   

Under this action, $\tau_x(u,g)$ is stable as $t \rightarrow \infty$ by proposition \ref{prop.F.stability}.  The stable subspace is trivial, so $\tau_x(u,g) = 0$ for almost-every $x \in M_0$, for all $(u,g) \in \nu_x^{(j),s}$.  
%The exponents of $(\Ad J^{(j-1)}_{\bf 0} F^t_x)^{-1}$ on $\mathfrak{x}^{(j-1)}$ are all at least $- \lambda^{(s)}$, so $\tau_x(u,g)$ must be trivial.  
%We conclude that $\mathfrak{y}_x^{(j)}(u,g) = \mathfrak{x}^{(j)}$.  
Induction on $j$ yields $\mathfrak{y}_x^{(r)}(u,g) = \mathfrak{x}^{(r)}$. 

\medskip
\emph{Step 4: Fibers of $\pi^r_0$ are connected, equal $\mathcal{X}^{(r)}$-orbits}.

For $x \in M_0$, write $B_x = B(\alpha(x)/2) \subset \mathcal{E}^{(r)}_x$.
From subsection \ref{subsubsec.digression.graph}, $\intnu_x^{(r),\lambda} \cap B_x$ is the graph of a function from a neighborhood of ${\bf 0}$ in $V_x^{(r),s}$ to $(V_x^{(r),s})^\perp$.  From proposition \ref{prop.Tr.spectrum} (2), $V_x^{(r),s} = \BR^n \oplus \widehat{V}^{(r),s}$.  Then $(\pi^r_0)^{-1}(u) \cap \intnu_x^{(r),\lambda} \cap B_x$ is the graph of a function on $\widehat{V}^{(r),s}$, so is connected.

%Let $\lambda$ and $\zeta$ be as in step 2, and choose $\lambda / 3 < \zeta < 0$. 
%% Given a compact set $K \subset (\pi^r_0)^{-1}(u) \cap \nu_x^{(r),s}$, there exists $k,m \geq 0$ such that
%% $$\interleave F^{(r)}(k+m,x) \cdot (u,g) \interleave \leq \beta(\varphi^k(x)) e^{m \lambda} \qquad \forall \ (u,g) \in K$$
Given $k \geq 0$, there exists $m \geq 0$ such that
$$\beta(\varphi^k(x)) e^{m \lambda} < \frac{C e^{ (k+m) \zeta}}{2} < \frac{\alpha(\varphi^{k+m}(x))}{2}$$ 
so 
$$ F^{(r)}(m,\varphi^k(x)) \cdot \intnu_{\varphi^k(x)}^{(r),\lambda} \subset B_{\varphi^{k+m}(x)}$$
Then
$$ \nu_x^{(r),s} \cap (\pi^r_0)^{-1}(u) = \bigcup_{k \geq 0} F^{(r)}(k,x)^{-1} \cdot \left((\pi^r_0)^{-1}(F^k_x u) \cap \intnu_{\varphi^k(x)}^{(r),\lambda} \cap  B_{\varphi^k(x)}  \right)$$
This is an increasing union of connected sets (see the proof of proposition \ref{prop.submanifolds.sequence}), so it is connected.  Now the fibers of $\nu_x^{(r),s}$ are connected $\mathcal{X}^{(r)}$-orbits, so $\nu_x^{(r),s}$ is a reduction to $\mathcal{X}^{(r)}$, invariant by $F^{(r)}(t,x)$.
\end{Pf}

The prolongation $F^{(r)}(t,x)$ is defined so that it preserves the section $({\bf 0}, \id)$.  The natural action of $F^t_x$ on $\mathcal{E}^{(r)}$ is by $(u,g) \mapsto (F^t_x (u), J^{(r)}_u F^t_x \cdot g)$.  Now we can interpret the $F^{(r)}(t,x)$-invariant submanifolds $\nu^{(r),s}_x$ in terms of this natural action.

\begin{corollary}
\label{cor.hr.reduction}
For any $r$ with $1 \leq r \leq q$, the natural action of $\{F^t_x \}$ on $\mathcal{E}^{(r)}$, and of any $(\{ G^k_x \},\psi) \in Z(F^t_x)$, preserves a reduction $\mathcal{R}^{(r)}$ to $\mathcal{H}^{(r)}$.   
\end{corollary}

\begin{Pf}
Let $\mathcal{R}^{(r)}$ be the saturation of $\cup_x \nu_x^{(r),s}$ by the right $\mathcal{H}^{(r)}$-action.  Recall from proposition \ref{prop.subres.subgroups} that $\mathcal{X}^{(r)} \lhd \mathcal{H}^{(r)}$.  By propositions \ref{prop.xr.reduction} and \ref{prop.J0.subres}, $\mathcal{R}^{(r)}$ is the desired $F^t_x$-invariant reduction.

For $(\{ G^k_x \}, \psi) \in Z(F^t_x)$, the prolongation $(\{ G^{(r)}(k,x) \}, \psi) \in Z(F^{(r)}(t,x))$.  Thus by proposition \ref{prop.submanifolds.centralizer}, the family $\cup_x \nu_x^{(r),s}$, restricted to an appropriate $\psi$-invariant, full measure subset, are invariant by $G^{(r)}(k,x)$.  By proposition \ref{prop.J0.subres}, $J_{\bf 0}^{(r)} G^k_x \in \mathcal{H}^{(r),0}$, so $\mathcal{R}^{(r)}$ is also $G^k_x$-invariant.
\end{Pf}

Fix $r = \lfloor \lambda^{(1)}/\lambda^{(s)} \rfloor$.
% there are no strict subresonances of $\Sigma^0$ of length greater than $r -1$, and no subresonances of length greater than $r$.  
We will next interpret the submanifolds $\cup_x \nu^{(r),s}_x$ as a $\varphi^t$-invariant family of rigid geometric structures on the fibers of $\mathcal{E}$.

\subsection{Dynamical submanifolds as geometric structures}

The reductions obtained in proposition \ref{prop.xr.reduction} and corollary \ref{cor.hr.reduction} above are not generally geometric structures of finite type in the sense of Cartan, because $\mathcal{X}^{(r)}$ contains rank-one elements (see \cite[Prop I.1.4]{kobayashi.transf}).  In this section, we implement a measurable version of Feres' approach in \cite{feres.normal.forms} to construct invariant generalized connections on the fibers $\mathcal{E}_x$ and coordinate atlases in which $F^t_x$ acts by resonance polynomials.  The automorphism group of each atlas on $\mathcal{E}_x$ is a finite-dimensional Lie group.  We further interpret these structures as an invariant family of flat connections on vector bundles over $\mathcal{E}_x$.

%% The reductions $\nu_x^{(r),s} \subset \mathcal{F}^{(r)}(\mathcal{E}_x)$ could be considered rigid geometric structures for the reason they carry a framing preserved by all $C^{r+1}$ local diffeomorphisms of $\mathcal{E}_x$ preserving $\nu_x^{(r),s}$, which makes these automorphisms into a Lie pseudogroup.  

Before stating the main theorem in general, we describe two important cases in detail.  Let $M_{\mathcal{R}}$ be the $\varphi^t$-invariant $\mu$-conull subset on which the conclusions of proposition \ref{prop.xr.reduction} hold.

\medskip

{\bf Example: $r=1$.}  This case arises when $\Sigma^0 = \{ \lambda \}$, or when $\lambda^{(s)} < \lambda^{(1)}/2$.  In the latter case the spectrum is said to be \emph{$1/2$-pinched}.  Then $\nu^{(1),s}_x$ is a $C^{q-1}$ reduction of $\mathcal{F}^{(1)}(\mathcal{E}_x)$ to $\{ \id \}$.  It determines a framing $g_x : \nu_x^s \rightarrow \GL(n)$ for all $x \in M_{\mathcal{R}}$.  Invariance by $F^{(1)}(t,x)$ means
$$ g_{\varphi^t(x)}(F^t_x(u)) = (J_u^{(1)} F^t_x) \cdot g_x(u) \cdot (J_{\bf 0}^{(1)} F^t_x)^{-1}; $$
in other words, 
$$ g_{\varphi^t(x)} \cdot T^t_x = (F^t_x)_* g_x$$
%% These framings $g_x$ on $\nu_x^s$ are thus $F^t_x$-invariant up to precomposition with $T^t_x \in \mathcal{H}^{(1),0} = \GL(W^1) \times \cdots \times \GL(W^s)$. 

The tangent spaces $T_\xi(\nu^{(1),s}_x)$ 
%is a horizontal space at $\xi \in \nu^{(1),s}_x$, 
project isomorphically to $T_{\pi^1_0(\xi)} \mathcal{E}_x$.  The union of right translates of these spaces by $\GL(n)$ is a $C^{q-2}$ horizontal distribution $\mathcal{D}_x$ on $\mathcal{F}^{(1)}(\mathcal{E}_x)$.  The $\oplus_x \mathcal{D}_x$ are principal connections, equivalent to $C^{q-2}$ affine connections $\oplus_x \nabla_x$ on $\cup_x \mathcal{E}_x$, invariant by $F^t_x$.
%which contains a neighborhood of ${\bf 0}$ in $\mathcal{E}_x$. 

%More precisely, $F(t,x)$ pushes forward $\nabla_x$ restricted to $\overline{B(\alpha(x))} \cap F(t,x)^{-1}(\overline{B(\alpha(\varphi^t(x)))}$ to $\nabla_{\varphi^t(x)}$ restricted to  $\overline{B(\alpha(\varphi^t(x)))} \cap F(t,x)(\overline{B(\alpha(x))})$; in particular, $F(t,\cdot)$ preserves the germs at ${\bf 0}$ of the family $\{ \nabla_x \}_{x \in M}$.   

The torsion $\oplus_x \tau_x$ of $\oplus_x \nabla_x$ is an $F^t_x$-invariant tensor. 
The spectrum of $T_{\mathcal{S}}$ on $\mathcal{S}_{\bf 0} = \oplus_x \wedge^2 T^*_{\bf 0}(\mathcal{E}_x) \otimes T_{\bf 0}(\mathcal{E}_x)$ is bounded below by $\lambda^{(1)} - 2 \lambda^{(s)} > 0$.  Propositions \ref{prop.perturbation.spectrum} and \ref{prop.F.stability} then imply that the connections $\nabla_x$ are torsion-free for all $x$ in a $\varphi^t$-invariant, $\mu$-conull subset of $M_{\mathcal{R}}$.

One could make a similar argument to show vanishing of the curvature of $\nabla_x$ assuming $q > 2$.  However, it is clear in any case that they are flat: the sections $\nu_x^{(1),s} \subset \mathcal{F}^{(1)}(\mathcal{E}_x)$ are parallel by construction.  Thus $\nu_x^{(1),s}$ correspond to coordinate charts $\theta_x$ on $\mathcal{E}_x$ for which $\theta_x^* \nabla_x$ equals the flat connection on $\BR^n$.  The resulting flat affine structures on $\cup_x \mathcal{E}_x$ are $F^t_x$-invariant; in fact, from the $F^t_x$-invariant framings above, we obtain a family of $\BR^n$-structures---framings by commuting vector fields---invariant by $F^t_x$ up to the linear action of the algebraic hull of $\{ T^t_x \}$ in $\mathcal{H}^{(1),0}$.

If $(\{G^k_x \}, \psi) \in Z(F^t_x)$, then $G^{(1)}(k,x)$ preserves $\cup_x \nu^{(1),s}_x$ and $\oplus_x \mathcal{D}_x$ (as usual, over $\cap_{i \in \BZ} \psi^i(M_{\mathcal{R}} \cap M_G)$), so the natural action of $\{G^k_x\}$ preserves the flat connections $\oplus_x \nabla_x$.  By proposition \ref{prop.J0.subres}, $\{ G^k_x \}$ also preserves the $\BR^n$-structures on $\cup_x \mathcal{E}_x$, up to the algebraic hull of $J_{\bf 0}^{(1)} G^k_x$, which is contained in $\mathcal{H}^{(1),0}$.

\begin{remark}
\label{rmk.parallel.filt.r=1}
We sketch how Ruelle's perturbation theorem implies that according to the frames in $\nu^{(1),s}_x$, the Lyapunov filtration is $V^1 \subset \cdots \subset V^s$ at every point of $\nu^s_x$: given $\xi \in \nu_x^{(1),s}(u)$, consider the sequences $\{ T_k \} = \{ T^1_{\varphi^{k-1}(x) } \}$ and $\{ T^{'}_k \} = \{ D_{F^{k-1}_x(u)} F^1_{\varphi^{k-1}(x)} \}$.  For $\lambda$ and $\eta$ as in section \ref{subsubsec.digression} 
%$\lambda' = \lambda + \eta < \lambda^{(i+1)}$, 
%appropriately chosen $\eta > 0$, 
%we have as in the proof of proposition \ref{prop.perturbation.spectrum} 
$$|| T_k' - T_k || e^{k \eta} < A e^{k \zeta}$$
 for some constant $A$ and $\lambda < \zeta < 0$, for all $k \geq 0$.  
%Let $W^{'1} \oplus \cdots \oplus W^{'m}$ be the Lyapunov decomposition 
Now \cite[(4.4)]{ruelle.diff.ergodic} implies that the $i$th Lyapunov projections $P^{i}$ for $\{ T^k \} = \{ T_k \circ \cdots \circ T_1 \}$ and $P^{'i}$ for $\{ T^{'k} \} = \{ D_uF^k_x \}$ satisfy
$$ || D_{u} F^\ell_x \circ P^{'i} \circ (D_{u} F^\ell_x)^{-1} - P^i || \leq A' e^{\ell \zeta} $$
for some constant $A'$, for all $\ell \geq 0$.  For $Q^{'i} = \xi^{-1} \circ P^{'i} \circ \xi$, 
%the $i$th Lyapunov projection in the frame at $u$ given by $\xi$,
$$ || D_{\bf 0} F^\ell_x \circ Q^{'i} \circ (D_{\bf 0} F^\ell_x)^{-1} - P^i || \leq B' e^{\ell \zeta} $$
for some $B'$, where we have used that $F^{(1)}(\ell,x) \cdot \xi$ tends to $\mbox{Id}$ at least as fast as $e^{\ell \lambda}$.
 It follows that the $i$th Lyapunov space at $u$ in the frame $\xi$ is congruent to $W^i$ modulo $V^i$ for all $i$, so the filtrations are equal.

%Denoting by $\mathcal{Q}_x$ 
Finally, we remark that $\mathcal{D}_x$ is tangent to the $\mathcal{H}^{(1),0}$-saturation of $\nu_x^{(1),s}$.
%,$\mathcal{Q}_x$.  
As $\mathcal{H}^{(1),0}$ preserves the Lyapunov filtration $V^1 \subset \cdots \subset V^s$, the Lyapunov filtration of $T \nu^{(1),s}_x$ is $\nabla_x$-parallel for all $x \in M_{\mathcal{R}}$.
\end{remark}
\medskip

{\bf Example: $r=2$.} In this case $3 \lambda^{(s)} < \lambda^{(1)} \leq 2 \lambda^{(s)}$.  
%We will again assume that the Lyapunov decomposition of $T_{\bf 0} \mathcal{E}$ is constant: $W^i_x \equiv W^i \subset \BR^n$ for $i = 1, \ldots, s$.

\subsubsection{Conglomeration of Lyapunov filtration to length $2$}  
\label{subsec.r=2.conglomeration}
Let $i_*$ be the minimal element of $\{1, \ldots , s\}$ such that
$\lambda^{(1)} - \lambda^{(i_*)} \leq \lambda^{(s)}$.  Whenever $j,l < i_*$ or $j,l \geq i_*$, then $\lambda^{(l)} - \lambda^{(j)} > \lambda^{(s)}$.  The first case is clear from the definition of $i_*$. For the second, note
$$ \lambda^{(l)} - \lambda^{(j)} \geq \lambda^{(i_*)}-\lambda^{(s)} \geq \lambda^{(1)}-2\lambda^{(s)} > \lambda^{(s)}$$

If we write $\widetilde{W}^1 = \oplus_{j=1}^{i_*-1} W^j$ and $\widetilde{W}^2 = \oplus_{j=i_*}^s W^j$, then 
$$ \mathcal{X}^{(1)} \subset \id + (\widetilde{W}^2)^* \otimes \widetilde{W}^1$$

% (** \emph{I can construct a similar coarsening of the Lyapunov filtration for larger $k$, but it doesn't say much about $\mathcal{X}_r$ for $r > 1$, so I'm not sure how useful it is}**)

\subsubsection{Framings of subquotients}
\label{subsec.framing.subquotient}

The $\mathcal{X}^{(1)}$-reductions $\nu_x^{(1),s} \subset \mathcal{E}^{(1)}_x$ induce distributions $\mathcal{V}^1_x$ on $\mathcal{E}_x$ corresponding to $\widetilde{W}^1$. For $\xi = (u,g) \in \nu_x^{(1),s}$ and $\xi' = F^{(1)}(t,x) \cdot \xi \in \nu_{\varphi^t(x)}^{(1),s}$,  
%Write $F^t_x \cdot \xi$ for the natural action of $F^t_x$ on $\mathcal{E}^{(1)}$; namely, $F^t_x \cdot \xi = (F^t_x (u), (J_u^{(1)} F^t_x) \cdot g)$ (we will use the same notation for the natural action of $F^t_x$ on $\mathcal{E}^{(r)}$).  Then
$$ \xi' \cdot T^t_x = F^t_x \cdot  \xi $$
%By proposition \ref{prop.J0.subres}, 
(the right-hand side is the natural action of $F^t_x$ on $\mathcal{E}^{(1)}$).
Since $T^t_x$ preserves $\widetilde{W}^1$ and $\widetilde{W}^2$, the distributions $\oplus_x \mathcal{V}^1_x$ are $F^t_x$-invariant.

Under restriction to $\widetilde{W}^1$ or projection modulo $\widetilde{W}^1$, 
$$ \mbox{Res}_{\widetilde{W}^1}(\mathcal{X}^{(1)}) = \{ \id \} \qquad \mbox{Proj}_{\BR^n/\widetilde{W}^1}(\mathcal{X}^{(1)}) = \{ \id \}$$
Therefore, $\cup_x \nu_x^{(1),s}$ also induces $F^{(1)}(t,x)$-invariant framings of $\oplus_x \mathcal{V}^1_x$ and $\oplus_x (T \mathcal{E}_x)/\mathcal{V}^1_x$.  If we denote $g_x$ the framing of $\mathcal{V}^1_x$, then
$$ (F^t_x)_* g_x = g_{\varphi^t(x)} \circ T^t_x$$
A similar identity holds for the framing of $T \nu_x^s / \mathcal{V}^1_x$.  
%Both framings are thus $F^t_x$-invariant, up to the linear actions given by $\mbox{Res}_{\widetilde{W}^1}(T^t_x)$ and $\mbox{Res}_{\widetilde{W}^2}(T^t_x)$. 

\subsubsection{Flat connections on $\cup_x \mathcal{E}_x$}
\label{subsec.r=2.connxns}

As $\mathcal{X}^{(2)} = \mathcal{X}^{(1)}$, the restriction $\pi^2_1: \nu_x^{(2),s} \rightarrow \nu_x^{(1),s}$ is a diffeomorphism.  The inverse induces a $C^{q-2}$ horizontal distribution $\mathcal{D}_x$ on $\nu_x^{(1),s}$ as follows.  Suppose that $\hat{\xi} = (u,J_{\bf 0}^{(2)} \varphi) \in \nu_x^{(2),s}$ projects to $\xi \in \nu_x^{(1),s}$.  Define a map $\BR^n \rightarrow \mathcal{E}^{(1)}_x$ 
$$ (\tau_u \circ \varphi , J^{(1)} \varphi) : v \mapsto  (\varphi(v)+ u, J_v^{(1)} \varphi)$$\
The subspace 
$$ S(\hat{\xi}) = \im D_{\bf 0}(\tau_u \circ \varphi, J^{(1)} \varphi) \subset T_\xi \mathcal{E}^{(1)}_x$$
depends only on $\hat{\xi}$.
Now set $\mathcal{D}_x(\xi) = S(\hat{\xi})$.  Note that this space projects under $(\pi^1_0)_*$ onto $T_{u} \mathcal{E}_x \cong \BR^n$.  Because $\mathcal{D}_x(\xi)$ comes from a diffeomorphism $\varphi$, it is called \emph{holonomic} (see \cite[p 29]{feres.framing.frobenius}).

For $h \in GL(n)$, write $R_{h*}$ for the derivative of the right action on $\mathcal{E}^{(1)}$.

\begin{proposition}
\label{prop.r=2.connxns}
The distribution $\mathcal{D}_x$ defines a $C^{q-2}$ principal connection $\nabla_x$ on $\nu_x^{(1),s}$ varying measurably in $x$; more precisely, for almost all $x \in M$, 
\begin{enumerate}
\item $\mathcal{D}_x(\xi) \subset T_\xi \nu_x^{(1),s}$
\item $(R_h)_* \mathcal{D}_x(\xi) = \mathcal{D}_x(\xi \cdot h)$ for all $h \in \mathcal{X}^{(1)}$.
\end{enumerate} 
The resulting connections $\nabla_x$ on $\mathcal{E}_x$ are flat.
\end{proposition}

\begin{Pf}
Write $\xi = (u,g)$ and $\hat{\xi} = (u, \hat{g}) = (u,J^{(2)}_{\bf 0} \varphi)$; set $\xi' = F^{(1)}(t,x) \cdot \xi$ and $\hat{\xi}' = F^{(2)}(t,x) \cdot \hat{\xi}$.  
%We will use $F^{(2)}(t,x)$-invariance of $\cup_x \nu_x^{(2),s}$.  
First we compare $\mathcal{D}_{\varphi^t(x)}(\xi') = S(\hat{\xi}')$ with $F^{(1)}(t,x)_* \cdot \mathcal{D}_x(\xi)$.  Write $(F^t_x)_*$ for the derivative of the natural action on $\mathcal{E}^{(1)}$ as in subsection \ref{subsec.framing.subquotient} (although it's ambiguous).  Set $\hat{h} = (J_{\bf 0}^{(2)} F^t_x)^{-1}$ and $h = \rho^2_1(\hat{h})$.  For the natural action, one can check
$$ S(F^t_x \cdot \hat{\xi}) = (F^t_x)_* S(\hat{\xi}) $$ 
For the principal right action, on the other hand, $S(\hat{\xi} \cdot \hat{h})$ is given by
\begin{equation}
\label{eqn.right.axn}
\begin{split}
D_{\bf 0}(\tau_u \circ \varphi \circ \hat{h}, J^{(1)}(\varphi \circ \hat{h})) & =  (R_h)_* \circ D_{\bf 0}(\tau_u \circ \varphi, J^{(1)} \varphi) \\
& + \omega^{-1} \circ D_{\bf 0} (J^{(1)} \hat{h})
\end{split}
\end{equation}
%%  is in the span of
%% $$ (F^t_x)_{* \xi} \circ R_{h*} (\mathcal{D}_{\hat{\xi}}) \qquad \mbox{and} \qquad  \omega_{F^{(1)}(t,x) \xi}^{-1}(\im (D_{\bf 0} J^{(1)} F^t_x))$$
(see \cite[eqn (54)]{feres.normal.forms}).  Because $\hat{h} \in \mathcal{H}^{(2),0}$ by proposition \ref{prop.J0.subres}, 
%% , the map
%% $$ (v,w) \mapsto \left( (D_{\bf 0} J^{(1)} F^t_x)(v) \right) (w)$$
%% corresponds to a degree two subresonance polynomial.
% which can only lie in $\mbox{Sym}^2 (\widetilde{W}^2)^* \otimes \widetilde{W}^1$.  
the image of $D_{\bf 0} (J^{(1)} \hat{h})$ must be strict subresonance---that is, in $\mathfrak{x}^{(1)}$.  In conclusion,
$$ S(\hat{\xi}') \equiv (F^t_x)_* \circ R_{h*} \cdot S(\hat{\xi}) \ \ \mbox{mod} \ \omega^{-1}(\mathfrak{x}^{(1)})$$
and thus
\begin{eqnarray}
\label{eqn.F1.on.D}
 \mathcal{D}_{\varphi^t(x)}(\xi') \equiv F^{(1)}(t,x)_* (\mathcal{D}_x(\xi)) \ \ \mbox{mod} \ \omega^{-1}(\mathfrak{x}^{(1)}).
\end{eqnarray}

Now define an $F^{(1)}(t,x)$-invariant family of sections $\tau_x$ of $(T \mathcal{E}_x)^* \otimes \left( \mbox{End} (T \mathcal{E}_x)/\omega^{-1}(\mathfrak{x}^{(1)}) \right)$, restricted to $\nu_x^{(1),s}$, by
$$ \mathcal{D}_x(\xi) \ \mbox{mod} \ \omega^{-1}(\mathfrak{x}^{(1)}) = \mbox{graph}(\tau_x({\xi}))$$

%where we identify $T \mathcal{E}_x \cong T \nu_x^{(1),s} / \omega^{-1}(\mathfrak{x}^{(1)})$
  The stability of $\tau$ given by proposition \ref{prop.F.stability} (together with proposition \ref{prop.perturbation.spectrum}) implies that $\tau_x$ vanishes for almost every $x$.  Thus $\mathcal{D}_x$ is tangent to $\nu_x^{(1),s}$ for almost every $x \in M$ and (1) is verified.

Point (2) follows from the fact that $\nu_x^{(2),s}$ and $\nu_x^{(1),s}$ are diffeomorphic reductions to $\mathcal{X}^{(1)} \subset \GL^{(1)}(n)$.  If $\hat{\xi}  \in \nu_x^{(2),s}$ projects to $\xi \in \nu_x^{(1),s}$, and $h \in \mathcal{X}^{(1)}$, then by equation (\ref{eqn.right.axn}),
$$ \mathcal{D}_x(\xi \cdot h) = S(\hat{\xi} \cdot h) = (R_h)_* \mathcal{D}_x(\xi).$$

Now that we have $\mathcal{D}_x$ tangent to $\nu_x^{(1),s}$ for almost every $x$, and given that $\nu_x^{(2),s}$ is at least $C^1$, we can argue as in \cite[sec 5.1]{feres.normal.forms}, \cite{feres.framing.frobenius} that the reductions $\nu_x^{(2),s}$ define a complete and consistent partial differential relation, and therefore, the distributions $\mathcal{D}_x$ are integrable.  It follows that the corresponding connections $\nabla_x$ on $\mathcal{E}_x$ are flat.  The integral leaves of $\mathcal{D}_x$ give a family of affine charts on $\mathcal{E}_x$ with transitions in the group $\mathcal{X}^{(1)}$.  
%We remark that since $\mathcal{X}^{(1)}$ preserves the filtration defined by $W^1 \oplus \cdots \oplus W^s$, it is parallel for $\nabla_x$.
\end{Pf}

\subsubsection{Flat connections on $\oplus_x \mathcal{V}^1_x$ and $\oplus_x (T \mathcal{E}_x)/\mathcal{V}^1_x$ invariant by $Z(F^t_x)$}

The connections $\oplus_x \nabla_x$ are not invariant by the natural action of $F^t_x$ on $\cup_x \mathcal{E}_x$.  By proposition \ref{cor.hr.reduction}, the natural action preserves an  $\mathcal{H}^{(2)}$-reduction $\mathcal{R}^{(2)} \subseteq \mathcal{E}^{(2)}$, containing the saturation of $\cup_x \nu_x^{(2),s}$ by the right action of $J_{\bf 0}^{(2)} F^t_x$.  We can extend $\mathcal{D}$ to a connection on $\mathcal{R}^{(1)} = \pi^2_1(\mathcal{R}^{(2)})$ using the right action of $\mathcal{H}^{(1)}$ (it is still not $F^t_x$-invariant).

We remark that there is an $F^t_x$-invariant horizontal distribution on $\mathcal{R}^{(2)}$.  It does not give a connection, but it is integrable.  See the general statements in theorem \ref{thm.main} (1) below.

From section \ref{subsec.framing.subquotient}, the distributions $\oplus_x \mathcal{V}^1_x$ are $F^t_x$-invariant.  There are thus $F^t_x$-equivariant maps $\mbox{Res}_{\mathcal{V}^1} : \mathcal{R}^{(1)} \rightarrow \mathcal{F} \mathcal{V}^1$ and $\mbox{Proj}_{T\mathcal{E}/\mathcal{V}^1} : \mathcal{R}^{(1)} \rightarrow \mathcal{F}(T \mathcal{E}/\mathcal{V}^1)$.  Note that $\omega^{-1}(\mathfrak{x}^{(1)})$ is in the kernel of both $(\mbox{Res}_{\mathcal{V}^1})_*$ and $(\mbox{Proj}_{T\mathcal{E}/\mathcal{V}^1})_*$.

Let $\xi \in \mathcal{R}^{(1)}$, with $F^t_x \cdot \xi = \xi'$.  Writing $h = J_{\bf 0}^{(1)} F^t_x$, we have
$$ (F^t_x)_* \mathcal{D}_x(\xi) = (R_h)_* \circ F^{(1)}(t,x)_*  \mathcal{D}_x(\xi)$$
 From equation (\ref{eqn.F1.on.D}) and the fact that $h$ normalizes $\mathcal{X}^{(1)}$,
$$ (F^t_x)_* \mathcal{D}_x(\xi) \ \mbox{mod } \omega^{-1}(\mathfrak{x}^{(1)}) \equiv (R_h)_* \mathcal{D}_{\varphi^t(x)}(\xi' \cdot h^{-1}) = \mathcal{D}_{\varphi^t(x)}(\xi')$$

Pushing forward $\mathcal{D}$ by $(\mbox{Red}_{\mathcal{V}^1})_*$ or $(\mbox{Proj}_{T\mathcal{E}/\mathcal{V}^1})_*$ thus gives $F^t_x$-invariant, flat connections on the vector bundles $\oplus_x \mathcal{V}^1_x$ and $\oplus_x T \mathcal{E}_x/V^1_x$.  Remark that by corollary \ref{cor.hr.reduction}, proposition \ref{prop.submanifolds.centralizer}, and proposition \ref{prop.J0.subres}, the distributions $\oplus_x \mathcal{V}_x^1$ and these flat connections are moreover invariant by all of $Z(F^t_x)$. Because $\mathcal{D}_x$ is tangent to $\mathcal{R}_x^{(1)}$ for almost all $x$ and $\mathcal{H}^{(1)}$ preserves the filtrations of $\widetilde{W}^1$ and $\BR^n / \widetilde{W}^1$ determined by $V^1 \subset \cdots \subset V^{i_*}$ and $V^{i_*}/\widetilde{W}^1 \subset \cdots \subset V^s/\widetilde{W}^1$, respectively, the Lyapunov filtrations of $\mathcal{V}^1_x$ and $T \mathcal{E}_x/\mathcal{V}^1_x$ are parallel for these connections, for almost all $x$ (see remark \ref{rmk.parallel.filt.r=1}).

%% Let the cocycle $F^t_x \in \mbox{Diff}^q(\BR^n, {\bf 0})$ over the ergodic flow $(\varphi^t,M,\mu)$ be as above, acting on the $\BR^n$-bundle $\mathcal{E} \rightarrow M$.  Assume  $ \sup_{0 \leq t \leq 1} \ln^+ \Vert F^{\pm t}_x \Vert_q \in L^1(M,\mu)$, where $q \geq 2$.  Assume that the Lyapunov spectrum of $T^t_x$ is 
%% $$ \lambda^{(1)} < \cdots < \lambda^{(s)} < 0$$
%% and set 

\bigskip

{\bf Main Theorem:} Let $r = \lfloor \lambda^{(1)}/\lambda^{(s)} \rfloor$ and $q \geq r+1$.

\begin{theorem}
\label{thm.main}
There exist the following differential-geometric structures on $\mathcal{E}_x$, for all $x$ in a $\varphi^t$-invariant, $\mu$-conull subset:
\begin{enumerate}
%% \item A $C^{q-r+1}$ reduction of $\mathcal{E}^{(r-1)}$ to $\mathcal{X}^{(r-1)}$, invariant by $F^t_x$, up to the algebraic hull of $J^{(r-1)}_{\bf 0} F^t_x$, which is contained in $\mathcal{H}^{(r-1),0}$.  If there are no resonances in $\Sigma^0$, then this algebraic hull is in $\GL(n)$.  This reduction is also invariant by any $(\{G^k_x \},\psi) \in Z(F^t_x)$, up to the algebraic hull of $J_{\bf 0}^{(r-1)} G^k_x$ in $\mathcal{H}^{(r-1),0}$.

\item A family $\mathcal{A}_x$ of $C^q$ charts, with transitions in $\mathcal{H}^{(r),0}$ (restricted to a neighborhood of ${\bf 0}$). The collection $\cup_x \mathcal{A}_x$ is $Z(F^t_x)$-invariant. 

\item A $C^{q-1}$ filtration $\mathcal{V}^1_x \subset \cdots \subset \mathcal{V}^l_x = T \mathcal{E}_x$, with $l \leq r$, equipped with $C^{q-2}$-smooth flat connections $\nabla_x^i$ on $\mathcal{V}_x^i/\mathcal{V}_x^{i-1}$, $i = 1, \ldots, l$.  The filtrations and connections are invariant by $Z(F^t_x)$. 
\end{enumerate}
\end{theorem}

\begin{Pf}
%% The reductions in (1) are given by proposition \ref{prop.xr.reduction}.  From theorem \ref{ruelle.submanifolds}, $\cup_x \nu_x^{(r-1),s}$ are invariant by the cocycle $F^{(r-1)}(t,x)$.  They are thus are invariant by the natural action of $F^t_x$ up to the right action of $J^{(r-1)}_{\bf 0} F^t_x$.  By proposition \ref{prop.J0.subres}, the latter has values in $\mathcal{H}^{(r-1),0}$.  As in the proof of corollary \ref{cor.hr.reduction}, these reductions are $G^{(r-1)}(k,x)$-invariant, and thus $G^k_x$-invariant up to the algebraic hull of $J^{(r-1)}_{\bf 0} G^k_x$, which, by proposition \ref{prop.J0.subres}, is also contained in $\mathcal{H}^{(r-1),0}$.
 Denote $\mathcal{R}^{(r)} = \left( \cup_x \nu_x^{(r),s} \right) \cdot \mathcal{H}^{(r)}$ the $Z(F^t_x)$-invariant $\mathcal{H}^{(r)}$-reduction of $\mathcal{E}^{(r)}$ given by corollary \ref{cor.hr.reduction}. We similarly have an invariant $\mathcal{R}^{(r+1)} \subset \mathcal{E}^{(r+1)}$, which is $C^{q-r-1}$-smooth.  As $\mathcal{H}^{(r)} = \mathcal{H}^{(r+1)}$, the restriction $\pi^{r+1}_r : \mathcal{R}^{(r+1)} \rightarrow \mathcal{R}^{(r)}$ is a diffeomorphism.  

We first construct a $Z(F^t_x)$-invariant family of holonomic distributions $\mathcal{D}_x(\xi) = S(\hat{\xi})$ on $\cup_x \nu_x^{(r),s}$ as in subsection \ref{subsec.r=2.connxns}; here $\hat{\xi}$ is the unique lift of $\xi \in \cup_x \nu_x^{(r),s} \subset \mathcal{R}^{(r)}$ to $\mathcal{R}^{(r+1)}$.  
%The construction of the atlases $\mathcal{A}_x$ is similar to the proof of part (1) of proposition \ref{prop.r=2.connxns}.  Namely, 
If $\hat{h} = J^{(r+1)}_{\bf 0} F^t_x \in \mathcal{H}^{(r+1),0}$, viewed as a local diffeomorphism fixing ${\bf 0}$, then $D_{\bf 0}(J^{(r)} \hat{h})$ has image in $\mathfrak{x}^{(r)}$.  We follow the proof of proposition \ref{prop.r=2.connxns} (1), and equation (\ref{eqn.F1.on.D}) becomes
$$ \mathcal{D}_{\varphi^t(x)}(\xi') \equiv F^{(r)}(t,x)_*(\mathcal{D}_x(\xi)) \ \mod \omega^{-1}(\mathfrak{x}^{(r)})$$
Then we can construct an invariant tensor $\sigma$ with values in the restriction to $\cup_x \nu_x^{(r),s}$ of
$$ \mathcal{T} = (T \mathcal{E}_x)^* \otimes ( \alggl^{(r)}(\mathcal{E}_x)/\omega^{-1}(\mathfrak{x}^{(r)}))$$
Here $\alggl^{(r)}(\mathcal{E}_x)$ denotes the tensor bundle over $\mathcal{E}_x^{(r)}$ with fiber 
$$\alggl^{(r)}_{(u,g)}(\mathcal{E}_x) = (\tau_u \circ g) \cdot \alggl^{(r)}(n) \cdot (\tau_u \circ g )^{-1}$$
Vanishing of $\sigma$ is implied by proposition \ref{prop.F.stability}, and we conclude that the restriction to $\nu_x^{(r),s}$ of $\mathcal{D}_x$ is tangent to $\nu_x^{(r),s}$ for all $x$ in a $\varphi^t$-invariant, $\mu$-conull set.

%The proof of part (2) of proposition \ref{prop.r=2.connxns} does \emph{not} carry over, because the right action of $\mathcal{X}^{(r)} = \mathcal{X}^{(r-1)}$ does not preserve $\mathcal{D}_x$ for $r > 2$.  
The distributions $\oplus_x \mathcal{D}_x$ are integrable, again by Gromov's Frobenius theorem, because $\cup_x \nu_x^{(r),s}$ is a complete and consistent partial differential relation; moreover, because they are holonomic, each leaf represents the $r$-jets of (the germ of) a diffeomorphism $(\BR^n,{\bf 0}) \rightarrow (\mathcal{E}_x, {\bf 0})$. 

Now we extend these distributions to $\mathcal{R}^{(r)}$ by the same formula: $\mathcal{D}_x(\xi) = S(\hat{\xi})$.  Any $\xi \in \mathcal{R}^{(r)}_x$ equals $\xi' \cdot h$ for $\xi' \in \nu_x^{(r),s}$ and $h \in \mathcal{H}^{(r)}$.  Equation (\ref{eqn.right.axn}) gives that $S(\hat{\xi}' \cdot \hat{h}) = \mathcal{D}_x(\xi)$ is tangent to $\mathcal{R}^{(r)}_x$.  Note that $\oplus_x \mathcal{D}_x$ are invariant by the natural $Z(F^t_x)$-action.  Given $\xi = ({\bf 0}, g) \in \mathcal{R}_x^{(r)}({\bf 0})$, denote $\alpha_\xi$ the germ of a coordinate parametrization of $\mathcal{E}_x$ with $\alpha_\xi({\bf 0}) = {\bf 0}$, $J_{\bf 0}^{(r)} \alpha_\xi = g$ and
$$ \im D_u (\alpha_\xi, J^{(r)} \alpha_\xi) = \mathcal{D}_x(\alpha_\xi(u),J_u^{(r)} \alpha_\xi) $$ 
for all $u$ in the domain of $\alpha_\xi$.  Note that this exists for all $\xi \in \nu_x^{(r),s}({\bf 0})$.

%Then set
% $$\mathcal{A}_x = \{ \alpha_\xi \ : \  \xi \in \mathcal{R}_x^{(r)}({\bf 0}) \}.$$
%Note that $\cup_x \mathcal{A}_x$ are $Z(F^t_x)$-invariant because $\cup_x \mathcal{R}_x^{(r)}$ and $\oplus_x \mathcal{D}_x$ are.

Given $\alpha_\xi$ as above for $\xi \in \mathcal{R}_x^{(r)}({\bf 0})$, for $\eta = (\alpha_\xi(u), J_u^{(r)} \alpha_\xi)$, the lift $\hat{\eta}$ to $\mathcal{R}^{(r+1)}$ equals $(\alpha_\xi(u), J_u^{(r+1)} \alpha_\xi)$: let $\varphi$ be a representative of $\hat{\eta}$.  Then $\eta = (\varphi({\bf 0}),J^{(r)}_{\bf 0} \varphi)$; in particular, $J_{\bf 0}^{(1)} \varphi = J_u^{(1)} \alpha_\xi$.  Also, $S(\hat{\eta})$ equals $\im D_{u}( \alpha_\xi , J^{(r)} \alpha_\xi)$.  The two equalities together imply $J_{\bf 0}^{(r+1)} \varphi = J_{u}^{(r+1)} \alpha_\xi$, which proves the claim.

% There are a couple ways to verify this; the easiest in our case is to note that our integrability argument works just as well on $\mathcal{R}^{(r+1)}$, which is locally diffeomorphic to a reduction $\mathcal{R}^{(r+2)}$ of $\mathcal{E}^{(r+2)}$ to $\mathcal{H}^{(r+2)} = \mathcal{H}^{(r+1)}$.  The holonomic integral leaves of the distributions $\hat{\mathcal{D}}_x$ on $\mathcal{R}^{(r+1)}$ project to the integral leaves of $\mathcal{D}_x$ in $\mathcal{R}^{(r)}$.  The coordinate parametrizations specified by these leaves are thus the same.  That means the image of $(\alpha_\xi, J^{(r+1)} \alpha_\xi)$ is a leaf in $\mathcal{R}^{(r+1)}$.

Now given $h \in \mathcal{H}^{(r)} = \mathcal{H}^{(r+1)}$ and $u$ in the domain of $\alpha_\xi \circ h$, let $\eta = (\alpha_\xi(h(u)), J_{h(u)}^{(r)} \alpha_\xi) \in \mathcal{R}_x^{(r)}$.  We can write
$$ ((\alpha_\xi \circ h)(u), J_u^{(r)} (\alpha_\xi \circ h)) = \eta \cdot (\tau_{- h(u)} \circ h \circ \tau_u)$$
Note that $\tau_{-h(u)} \circ h \circ \tau_u \in \mathcal{H}^{(r+1)}$.  The lift of our point to $\mathcal{R}^{(r+1)}$ is 
$$  (\alpha_\xi (h(u)), J_{h(u)}^{(r+1)} \alpha_\xi)  \cdot (\tau_{- h(u)} \circ h \circ \tau_u) = ((\alpha_\xi \circ h)(u), J_u^{(r+1)} (\alpha_\xi \circ h))$$
The $r$-jets of $\alpha_\xi \circ h$ thus lie in the integral leaf of $\mathcal{D}_x$ through $\xi \cdot h$.
% $$ D_{\bf 0} (\alpha_\xi \circ \tau_u \circ h, J^{(r)} (\alpha_\xi \circ \tau_u \circ h)) = S(\hat{\eta} \cdot h) = \mathcal{D}_x(\eta \cdot h)$$
We conclude $\alpha_{\xi h} = \alpha_\xi \circ h$ for all $\xi \in \mathcal{R}^{(r)}_x({\bf 0})$ and $h \in \mathcal{H}^{(r)}$.  In particular, given $\alpha_\xi$ for $\xi \in \nu_x^{(r),s}({\bf 0})$, we can define $\alpha_{\xi \cdot h} = \alpha_\xi \circ h$.

Now $Z(F^t_x)$ preserves $\cup_x \mathcal{A}_x$, and each $\mathcal{A}_x$ is parametrized by $\mathcal{R}^{(r)}_x({\bf 0}) = \{ {\bf 0} \} \times \mathcal{H}^{(r)}$, so the action of $F^t_x$ on $\mathcal{E}_x$ is determined by $J^{(r)}_{\bf 0} F^t_x$ (the same is true for any $\{ G^k_x \}$ in $Z(F^t_x)$).  By proposition \ref{prop.J0.subres}, the subatlas
$$ \mathcal{A}^0_x = \{ \alpha_\xi \ : \ \xi \in \{ {\bf 0} \} \times \mathcal{H}^{(r),0} \}$$

is $Z(F^t_x)$-invariant.  The transitions between the charts in $\mathcal{A}^0_x$ belong to $\mathcal{H}^{(r),0}$.  Point (1) is proved.

The first step to prove part (2) is to conglomerate the Lyapunov filtration as in subsection \ref{subsec.r=2.conglomeration}.  Recall that $(r+1) \lambda^{(s)} < \lambda^{(1)} \leq r \lambda^{(s)}$.  Set $i_0 = 1$, and recursively define $i_k$ to be the minimal $i > i_{k-1}$ with $\lambda^{(i_{k-1}) } - \lambda^{(i_k)} \leq \lambda^{(s)}$.  
%Let $r_* \leq r$ be such that the process terminates with $i_{r_*-1}$.  
It is a simple induction argument to see that $\lambda^{(i_k)} - (r-k+1) \lambda^{(s)} > 0$.  It follows that the process terminates with $i_{r_*-1}$ where $r_* -1 < r$.  Set $i_{r*} = s+1$, and 
$$ \widetilde{W}^l = \bigoplus_{j=i_{l-1}}^{i_l-1} W^j$$

Then $\widetilde{W}^1 \oplus \cdots \oplus \widetilde{W}^{r_*}$ is an $\mathcal{H}^{(1),0}$-invariant decomposition, of length at most $r$, with respect to which $\mathcal{X}^{(1)}$ is block upper-triangular.  Then, as in section \ref{subsec.framing.subquotient}, we obtain $C^{q-1}$ filtrations
$$ 0 \subset \mathcal{V}_x^1 \subset \cdots \subset \mathcal{V}_x^{r_*} = T \mathcal{E}_x$$
invariant by $Z(F^t_x)$.  For $1 \leq l \leq r$, define a map on $\mathcal{R}^{(r)}$
$$ Q_l = \mbox{Proj}_{\mathcal{V}^l/\mathcal{V}^{l-1}} \circ \mbox{Res}_{\mathcal{V}^l} \circ \pi^r_1$$

Note that $\mathcal{D}_x$ is not projectible under $(\pi^r_1)_*$ because it is not right $\mathcal{H}^{(r)}$-invariant; however, it does descend under $(Q_l)_*$ to a horizontal distribution on $\mathcal{F}(\mathcal{V}^l/\mathcal{V}^{l-1})$: The degree $r$ analogue of equation (\ref{eqn.right.axn}) says that $\mathcal{D}_x$ is right-$\mathcal{H}^{(r)}$-invariant modulo $\omega^{-1}(\mathfrak{x}^{(r)})$.  But $\omega^{-1}(\mathfrak{x}^{(1)})$ is in the kernel of $(\mbox{Proj}_{\mathcal{V}^l/\mathcal{V}^{l-1}} \circ \mbox{Res}_{\mathcal{V}^l})_*$.  By the same reasoning, the projection of $\mathcal{D}_x$ to a horizontal distribution on $\mathcal{F}(\mathcal{V}^l/\mathcal{V}^{l-1})$ gives a $C^{q-2}$ principal connection $\overline{\nabla}^l$.  The principal group here is $\GL(\widetilde{W}^l)$, which can be obtained as a quotient of $\mathcal{H}^{(r)}$.  Because $\mathcal{D}_x$ is integrable, so is the projection by $(Q_l)_*$; thus $\overline{\nabla}^l$ is flat.

As $\mathcal{D}_x$ is tangent to $\mathcal{R}_x$, the Lyapunov filtration of $\mathcal{V}_x^\ell / \mathcal{V}_x^{\ell -1}$ is $\overline{\nabla}^\ell$-parallel for all $\ell$, for almost all $x$, as in remark \ref{rmk.parallel.filt.r=1} and subsection \ref{subsec.r=2.connxns}.
\end{Pf}

\subsection{Foliated case: Smooth geometric structures on leaves}

We return to the important special case that $M$ is a compact $C^0$ manifold and the fibers $\mathcal{E}_x$ are plaques of a $\varphi^t$-invariant foliation; we assume that the action on this foliation is contracting.   More precisely, let $L$, as in section \ref{sec.filtered.foliations}, be a $C^0$ foliation of $M$ by $n$-dimensional submanifolds admitting a $\varphi^t$-invariant $C^q$ smooth structure, and assume $\sup_{0 \leq t \leq 1} \interleave J_x^{q} (\left. \varphi^{\pm t} \right|_L) \interleave$ is bounded in $x$ (for example, it is continuous in $x$). 
%\ln^+ \Vert \left. \varphi^{\pm t} \right|_{L} \Vert_{q,x} \in L^1(M,\mu)$.  
Let $\{ (\mathcal{E}_x, \theta_x)\}$ be a uniformly biLipschitz $C^q$ atlas along $L$.  Define the cocycle $\{ F^t_x \}$ as in section \ref{sec.filtered.foliations}, and assume that the Lyapunov exponents for $\{ T^t_x = D_{\bf 0} F^t_x \}$ are all negative. 

Before interpreting the results of the previous section in this setting, we must assemble the manifolds $\nu_x^{(r),s}$ associated to the prolongation of $\{F^t_x\}$ into an invariant submanifold of the $r$-frame bundle $\cup _x \mathcal{F}^{(r)} L_x = \cup_x L^{(r)}_x$ along $L$.  

\begin{proposition}
\label{prop.prolonged.foliations}
Let $r \leq q-1$. For all $x$ in a $\varphi^t$-invariant, $\mu$-conull subset of $M$, 
\begin{enumerate}
\item there exists a neighborhood $x \in U_x \subset L_x$, and a $C^{q-r}$-smooth reduction $\mu_x^{(r),s}$ of $\left. L^{(r)}_x \right|_{U_x}$ to $\mathcal{X}^{(r)}$, such that
\item $\cup_{y \in L_x} \mu_y^{(r),s} \cdot \mathcal{H}^{(r)}$ is a $C^{q-r}$-smooth $\mathcal{H}^{(r)}$-reduction $\mathcal{P}^{(r)}_x \subset L_x^{(r)}$, and $\cup_x \mathcal{P}^{(r)}_x$ is invariant by $Z(\varphi^t)$. 
\end{enumerate}
\end{proposition}

  %% We will start with an atlas along $L$ that is appropriately bounded and use this to assemble submanifolds in a bounded open subset of $\cup_x L^{(r)}_x$.   Recall that for proposition \ref{prop.xr.reduction}, we applied measurable linear changes of coordinates on $\mathcal{E}$ in order to make the Lyapunov decomposition constant.  These changes may not be bounded over $M$.  We will nonetheless perform these fiberwise linear changes and show that the resulting submanifolds still fit nicely together and retain their dynamical properties.   

%% The charts $\{ (\mathcal{E}_x^{(r)}, \theta^{(r)}_x) \}$ will not in general give a cocycle $T^{(r)}$ with constant Lyapunov decomposition.  For $y \in L_x^s$, the point $\theta^{(r)}_y({\bf 0}, \id)$ cannot be expected to be in $\theta^{(r)}_x (\nu_x^{(r),s})$, because these two charts may have nothing to do with each other. 

\begin{Pf} \emph{Step 1: Bounded atlas along $L^{(r)}$.}

Some extension of the uniformly biLipschitz property is required, although the fibers of $L^{(r)}$ are not compact.  We choose a $C^q$ atlas $\{ (\mathcal{E}_x, \tilde{\theta}_x) \}$ along $L$ with the properties:
\begin{itemize}
\item  The prolongations $\tilde{\theta}^{(r)}_x : \mathcal{E}^{(r)}_x \rightarrow L^{(r)}_x$ are such that $\cup_x \tilde{\theta}_x^{(r)}(\overline{B(1)})$ lies in a compact subset $K$ of $\cup_x L_x^{(r)}$;
\item  $\{(\mathcal{E}_x^{(r)}, \tilde{\theta}_x^{(r)}) \}$ is uniformly biLipschitz with respect to the Finsler metric $d$ on $\mathcal{E}^{(r)}$ and some continuous metric $d_K$ on $K$.
\end{itemize}
(As in section \ref{sec.filtered.foliations}, such an atlas can be obtained from an appropriately chosen finite cover of $M$ by foliated charts.)
Note that the resulting $C^{q-r}$ atlas along $L^{(r)}$ is $\GL^{(r)}(n)$-equivariant.  
%View $L^{(r)}$ as a $C^0$ foliation of the bundle $\cup_x L_x^{(r)} \rightarrow M$.  
Let $\{\widetilde{F}^t_x\}$ be the cocycle determined by $\{ (\mathcal{E}_x, \tilde{\theta}_x) \}$.  The prolonged cocycle satisfies (compare sections \ref{sec.filtered.foliations}, \ref{subsec.xr.reduction})
$$ \sup_{0 \leq \epsilon \leq 1} \ln^+ \Vert \widetilde{F}^{(r)}(\pm \epsilon,x) \Vert_{q-r} \in L^1(M,\mu).$$
The set $\{ g : Kg \cap K \neq \emptyset \}$ lies in a compact $C \subset \GL^{(r)}(n)$.
% there is also a uniform bound on the Lipschitz constants of elements of $C$ restricted to $K$.  
Let $m$ be a Lipschitz constant valid for the atlas and for all $g \in C$ restricted to $K$.  

%% It is also uniformly biLipschitz. (As in section \ref{sec.filtered.foliations}, $\omega$ defines a left-invariant Riemannian metric $d_G$ on the fibers of $\mathcal{E}^{(r)} \cong \BR^n \times \GL^{(r)}(n)$; ... AND??
% the following metric on $L^{(r)}_x$ is thus well defined:
% $$ \hat{d}(\theta^{(r)}_x (u,g),\theta^{(r)}_x(v,h)) = \left( d(\theta_x(u),\theta_x(v))^2 + d_G(g,h)^2 \right)^{1/2}$$ 
% and $\theta^{(r)}_x$ is isometric on fibers.)

\medskip
\emph{Step 2: Construction of $\tilde{\mu}_x^{(r),s}$ from bounded atlas.}
 
Let $\tilde{\nu}_x^{(r),s}$ be the submanifolds given by theorem \ref{thm.ruelle.submanifolds} applied to $\{ \widetilde{F}^{(r)}(t,x)\}$. Write $B = \overline{B(1)}$.
Let $(u,g) \in \tilde{\nu}_x^{(r),s} \cap B$ with $\tilde{\theta}_x(u) = y \in L_x$, and
$$\tilde{\theta}^{(r)}_x(u,g) = \eta = \tilde{\theta}^{(r)}_y({\bf 0}, \id) \cdot h \qquad h \in C \subset \GL^{(r)}(n)$$ 

  We wish to show that 
\begin{equation}
\label{eqn.nux.in.nuy}
\tilde{\theta}^{(r)}_x \left( \tilde{\nu}_x^{(r),s} \cap B \right) \ \cap \ \tilde{\theta}_y^{(r)} \left( B \cdot h \right) \subset \tilde{\theta}_y^{(r)} \left( \tilde{\nu}_y^{(r),s} \cdot h \right)
\end{equation}
Let  $\xi = \tilde{\theta}_x^{(r)} (u',g')  =  \tilde{\theta}^{(r)}_y(v,h')$ with $(u',g') \in \tilde{\nu}_x^{(r),s} \cap B$ and $(v,h'h^{-1}) \in B$.  By the biLipschitz property of $\tilde{\theta}^{(r)}_{\varphi^t(y)}$,
\begin{multline}
\label{eqn.step2} 
d\left( (\widetilde{F}^t_y(v), (J_{v}^{(r)} \widetilde{F}^t_y) \cdot h'h^{-1} \cdot (J_{\bf 0}^{(r)} \widetilde{F}^t_y)^{-1}),({\bf 0}, \id) \right) \\ \leq m d_K \left( \varphi^t \xi \cdot h^{-1} (J_{\bf 0}^{(r)} \widetilde{F}^t_y)^{-1}, J_{\bf 0}^{(r)} \tilde{\theta}_{\varphi^t(y)} \right) 
\end{multline}

The following points are both in $K$ (for sufficiently large $t > 0$):
$$\varphi^t \eta \cdot h^{-1} (J_{\bf 0}^{(r)} \widetilde{F}^t_y)^{-1} = J_{\bf 0}^{(r)} \tilde{\theta}_{\varphi^t(y)} \ \ ; \ \ \varphi^t \eta \cdot (J_{\bf 0}^{(r)} \widetilde{F}^t_x)^{-1} = \tilde{\theta}^{(r)}_{\varphi^t x} \left( \widetilde{F}^{(r)}(t,x) (u,g) \right) $$

Now right translate by $(J_{\bf 0} \widetilde{F}^t_y) \cdot h \cdot (J_{\bf 0} \widetilde{F}^t_x)^{-1} \in C$ 
to obtain the bound
$$ (\ref{eqn.step2}) \leq m^2 d_K \left( \varphi^t \xi \cdot (J_{\bf 0}^{(r)} \widetilde{F}^t_x)^{-1}, \varphi^t \eta \cdot (J_{\bf 0}^{(r)} \widetilde{F}^t_x)^{-1} \right), $$
and apply $(\tilde{\theta}_{\varphi^t(x)}^{(r)})^{-1}$ to obtain
$$\leq m^3 d \left( \widetilde{F}^{(r)}(t,x)(u',g') , \widetilde{F}^{(r)}(t,x)(u,g) \right)  $$
Thus finally
$$\limsup_{t \rightarrow \infty} \frac{1}{t} \ln d \left( \widetilde{F}^{(r)}(t,y)(v,h'h^{-1}), ({\bf 0}, \id) \right) \leq \lambda^{(s)}$$

%%  We have 
%% $$ \lim \frac{1}{t} \ln \Vert (F^t_x(u), (J_{u}^{(r)} F^t_x) \cdot g \cdot (J_{\bf 0}^{(r)} F^t_x)^{-1}) \Vert \leq \lambda^{(s)}$$

and $(v,h'h^{-1}) \in \tilde{\nu}_y^{(r),s}$.  A similar argument shows 
\begin{equation}
\label{eqn.nuy.in.nux}
\tilde{\theta}^{(r)}_y \left( \tilde{\nu}_y^{(r),s} \ \cap \ B \right) \cdot h  \ \cap \ \tilde{\theta}_x^{(r)} \left( B \right) \subset \tilde{\theta}_x^{(r)} \left( \tilde{\nu}_x^{(r),s} \right)
\end{equation}
In conclusion, the manifolds $\tilde{\nu}_x^{(r),s} \cap B$ map forward under $\tilde{\theta}^{(r)}_x$ to submanifolds $\tilde{\mu}_x^{(r),s} \subset L^{(r)}_x$ which smoothly fit together with $\tilde{\mu}_y^{(r),s}$, for $y \in L_x$, after a vertical translation, according to (\ref{eqn.nux.in.nuy}) and (\ref{eqn.nuy.in.nux}). 

\medskip
\emph{Step 3: Tempered cocycle equivalence, construction of $\mu_x^{(r),s}$}

The Oseledec-Pesin reduction theorem \cite[Thm 6.10]{barreira.pesin.introset} gives a linear tempered equivalence $\{ g_x \}$ of $\{ \widetilde{T}^t_x \}$ with a cocycle $\{ T^t_x \}$ having constant Lyapunov decomposition (as we assumed in section \ref{subsec.Tr.spectrum}).  By proposition \ref{prop.Tr.spectrum}, the Lyapunov decomposition of $\{ T^{(r)}(t,x) \}$ is also constant, for $x \in M_T$, defined as in section \ref{subsec.Tr.spectrum}.  
Recall that tempered means 
$$ \lim_{t \rightarrow \pm \infty} \frac{1}{t} \ln \Vert g_{\varphi^t(x)}^{\pm 1} \Vert  = 0$$

% Let $M_0$ as in theorem \ref{thm.ruelle.submanifolds}, where $M_F$ and $M_T$ are defined in terms of $\{ \widetilde{F}^{(r)}(t,x) \}$ and $\{ T^{(r)}(t,x) \}$, respectively.  Fix a decomposition $\BR^n = W^1 \oplus \cdots \oplus W^s$ isomorphic to the Lyapunov decomposition for $T$ in $M_0$, and let $g_x (W^i_x) = W^i$ for $i=1,\ldots, s$, where $g_x \in GL(n)$ varies measurably in $x$.  By proposition \ref{prop.Tr.spectrum}, $g_x$ transforms the Lyapunov decomposition of $T^{(r)}$ at $x \in M_0$ to $\oplus_{\sigma} (W^\sigma \oplus \widehat{W}^{(r),\sigma})$.  The norm $\Vert g_x \Vert$ can be bounded above in terms of the angles between the Lyapunov spaces at $x$.  By REF, these tend to zero at a subexponential rate, which implies 

Now set $\theta_x = \tilde{\theta}_x \circ g_x$.  The new cocycles are $F^t_x = g_{\varphi^t(x)}^{-1} \circ \widetilde{F}^t_x \circ g_x$ and 
$$F^{(r)}(t,x)(u,h) = \left( g_{\varphi^t(x)}^{-1} \circ \widetilde{F}^{(r)}(t,x)(g_x u,g_x h g_x^{-1}) \right) \cdot g_{\varphi^t(x)}$$

By proposition \ref{prop.poly.bounds} and the tempered property,
\begin{eqnarray*}
& & \limsup_{t \rightarrow \infty} \frac{1}{t} \ln d \left( \widetilde{F}^{(r)}(t,x)(g_x u,g_x h g_x^{-1}), ({\bf 0}, \id) \right)  \leq \lambda^{(s)}  \ \Leftrightarrow \\
& &  \limsup_{t \rightarrow \infty} \frac{1}{t} \ln d \left( F^{(r)}(t,x)(u,h), ({\bf 0}, \id) \right) \leq \lambda^{(s)}
\end{eqnarray*}
Then 
$$\nu_x^{(r),s} = (g_x^{-1}, \Ad g_x^{-1}) \cdot \tilde{\nu}_x^{(r),s} \qquad \mbox{and} \qquad \theta_x^{(r)} \left( \nu_x^{(r),s} \right) = \tilde{\mu}_x^{(r),s} \cdot g_x$$
By proposition \ref{prop.xr.reduction}, the submanifolds $\mu_x^{(r),s} = \tilde{\mu}_x^{(r),s} \cdot g_x$ are $C^{q-r}$-smooth $\mathcal{X}^{(r)}$-reductions of $L_x^{(r)}$ over $U_x = \theta_x(\mathcal{E}_x)$, for all $x$ in a $\varphi^t$-invariant, $\mu$-conull set $M_{\mathcal{R}}$.  They smoothly fit together when appropriately vertically translated.  Point (1) is proved.

\medskip
\emph{Step 4: $\mathcal{H}^{(r)}$ saturation, construction of $\mathcal{P}_x^{(r)}$}

Let $y \in \theta_x(\mathcal{E}_x)$.  From above, there exists $g \in \GL^{(r)}(n)$ such that $\mu_x^{(r),s} \cap \mu_y^{(r),s} \cdot g$ is open in each term.  A point $\xi \in  \mu^{(r),s}_x$ determines a frame of $\mu_x^{(r-1),s}$ at $\overline{\xi} = \pi^r_{r-1} (\xi)$ in which, as in remark \ref{rmk.parallel.filt.r=1}, the Lyapunov filtration comprises the subspaces $V^{(r-1),\sigma}$ for $\sigma \leq \lambda^{(s)}$.  On the other hand, the arguments of step 2 above show that the derivative of right translation by $\overline{g} = \rho^r_{r-1}(g)$ carries the Lyapunov filtration of $\mu_y^{(r-1),s}$ at $\overline{\xi}$ to that of $\mu_x^{(r-1),s}$ at $\overline{\xi} \cdot \overline{g}$.  It follows (see formula (\ref{eqn.right.axn}) and its higher-order analogues) that
%The Lyapunov filtrations of $T\mu_x^{(r-1),s}$ and $T\mu_y^{(r-1),s}$ are related by $g$, which implies 
$g \in \mathcal{H}^{(r)}$.  For $x \in M_{\mathcal{R}}$, set $\mathcal{P}_x^{(r)} = \cup_{y \in L_x} \mu_y^{(r),s} \cdot \mathcal{H}^{(r)}$, a $C^{q-r}$ reduction of $L_x^{(r)}$ to $\mathcal{H}^{(r)}$.
 
Invariance of $\cup_x \nu_x^{(r),s}$ by $F^{(r)}(t,x)$ gives
$$ \varphi^t \mu_x^{(r),s} = \mu_{\varphi^t(x)}^{(r),s} \cdot (J_{\bf 0}^{(r)} F^t_x)$$
By proposition \ref{prop.J0.subres}, $J_{\bf 0}^{(r)} F^t_x \in \mathcal{H}^{(r),0}$, so $\cup_x \mathcal{P}^{(r)}_x$ are $\varphi^t$-invariant.  Any $\psi \in Z(\varphi^t)$ (see remark \ref{rmk.foliations.centralizer}) enjoys a similar invariance of the $\mu_x^{(r),s}$, and by proposition \ref{prop.J0.subres}, preserves $\cup_x \mathcal{P}^{(r)}_x$ where $x \in \cap_{i \in \BZ} \psi^i(M_{\mathcal{R}} \cap M_G)$. \end{Pf}

Our first theorem for contracted foliations concerns the atlas from theorem \ref{thm.main} (1): it gives rise to a homogeneous structure on leaves. 

\begin{theorem}
\label{thm.foliation.atlas}
Assume $q \geq \lfloor \lambda^{(1)}/\lambda^{(s)} \rfloor + 1$.  For all $x \in M_{\mathcal{R}}$, for all $y \in L_x$, there is a family  $\mathcal{B}^0(y)$ of $C^q$ charts on $L_x$, varying measurably between leaves and with the following additional properties:
% There is a $Z(\varphi^t)$-invariant, measurable family of atlases $\cup_x \mathcal{B}_x^0$ on the leaves of $L$, where $\mathcal{B}_x^0 = \cup_{y \in L_x} \mathcal{B}_x^0(y)$, such that:
\begin{enumerate} 
\item Global diffeomorphisms: For all $x$ in a $\varphi^t$-invariant, $\mu$-conull set, each $\beta \in \mathcal{B}^0(y)$ is a diffeomorphism $\beta : (\BR^n, {\bf 0}) \rightarrow (L_x,y)$.
\item Homogeneity: Given $z \in L_x$, $\beta \in \mathcal{B}^0(y)$, and $\gamma \in \mathcal{B}^0(z)$, we have $(\gamma \circ \beta^{-1}) \cdot \mathcal{B}^0(y) = \mathcal{B}^0(z)$.
\item Invariance: Let $\mathcal{B}^0_x = \cup_{y \in L_x} \mathcal{B}^0(y)$.  The collection $\cup_x \mathcal{B}_x^0$ is $Z(\varphi^t)$-invariant, and in these coordinates, $Z(\varphi^t)$ acts by resonance polynomials---that is, elements of $\mathcal{H}^{(r),0}$.
\item Structure group: For $\beta, \gamma \in \mathcal{B}_x^0$, the transition $\gamma^{-1} \circ \beta$ acts on $\BR^n$ by a translation composed with the an element of $\mathcal{H}^{(r)}$; in particular, $\mathcal{B}_x^0$ forms a $C^q$ atlas on $L_x$.
\end{enumerate}
If the spectrum $\Sigma^0$ is $1/2$-pinched, then $L$ carries an invariant family of $C^{q-1}$ flat affine structures, measurable in $x$.
\end{theorem}

As usual, the $Z(\varphi^t)$-invariance in (3) holds for a given centralizing $\psi$ on an appropriate conull, $\psi$-invariant subset of $M_{\mathcal{R}}$.

\begin{Pf}
Set $r = \lfloor \lambda^{(1)}/\lambda^{(s)} \rfloor$, and let $\cup_x \mathcal{P}_x^{(r)} = \cup_x \mu_x^{(r),s} \cdot \mathcal{H}^{(r)}$ be the reduction of $L^{(r)}$ to $\mathcal{H}^{(r)} \cong \mathcal{H}^{(r+1)}$ given by proposition \ref{prop.prolonged.foliations} (2).  Let $\{ (\mathcal{E}_x, \theta_x) \}$ be the foliated atlas along $L$ from the proof.

As in the proof of theorem \ref{thm.main} (1), define integrable horizontal distributions on $\mu_x^{(r),s}$ and extent to integrable, $Z(\varphi^t)$-invariant, horizontal distributions on $\mathcal{P}_x^{(r)}$.  Given $\xi \in \mathcal{P}^{(r)}_x$, denote $\widehat{L}_\xi$ the integral leaf through $\xi$.  For $\xi = \theta^{(r)}_y ({\bf 0},g)$, let $\beta_\xi = \theta_y \circ \alpha_{({\bf 0},g)}$, where $\alpha_{({\bf 0},g)}$ is as in the proof of theorem \ref{thm.main} (1).  Note that the image of $J^{(r)} \beta_\xi$ is in $\widehat{L}_\xi$.  Set
$$\mathcal{B}^0(y) = \{ \beta_\xi \ : \ \xi \in \theta^{(r)}_y(\{ {\bf 0} \} \times \mathcal{H}^{(r),0}) \} $$

%For $y \in L_x$, the fiber $(\theta^{(r)}_y)^{-1} \mathcal{P}^{(r)}_x(y) = \{ {\bf 0} \} \times \mathcal{H}^{(r)}$.

For $x \in M_{\mathcal{R}}$, write $\hat{x} = J_{\bf 0}^{(r)} \theta_x$.  The $C^{q-r}$ leaves $\widehat{L}_{\widehat{x}} \subset L^{(r)}_x$ vary measurably in $x$, and therefore so does the maximum domain of definition of $\beta_{\hat{x}}$.  Lusin's theorem gives $\epsilon > 0$ and a compact $C \subset M$, visited infinitely many times by $\{ \varphi^t(x) \}$ for almost all $x$, and such that $\beta_{\hat{y}}$ is defined on $B(\epsilon)$ for all $y \in C$.

Let $J_{\bf 0}^{(r)} F^t_x = h_t \in \mathcal{H}^{(r),0}$, which are diffeomorphisms by proposition \ref{prop.subres.subgroups}.  By remark \ref{rmk.rjets.contract}, $h_k \rightarrow {\bf 0}$ uniformly on compact sets as $k \rightarrow \infty$ in $\BN$.  Next,
$$ \beta_{\hat{x}} = \varphi^{-k} \circ \beta_{\widehat{\varphi^k(x)}} \circ h_k$$
Let $x \in M_{\mathcal{R}}$ and $k_m \rightarrow \infty$ be such that $\varphi^{k_m}(x) \in C$ for all $m \geq 1$.  The domain of $\beta_{\hat{x}}$ contains $\cup_m h_{k_m}^{-1}( B(\epsilon)) = \BR^n$.  Given $y \in L_x$, there is $m$ such that $\varphi^{k_m}(y) \in \im \beta_{\widehat{\varphi^{k_m}(x)}}$.  Thus $y \in \im \beta_{\hat{x}}$.  Now $\beta_{\hat{x}}$ is a diffeomorphism as claimed in (1).  As in the proof of theorem \ref{thm.main} (1), $\beta_{\xi \cdot h} = \beta_\xi \circ h$, so $\beta_\xi$ is a diffeomorphism for all $\xi$ in $\mathcal{P}_x^{(r)}$, for all $x \in M_{\mathcal{R}}$.

Let $\beta = \beta_\xi \in \mathcal{B}_x^0(y), \gamma = \beta_\omega \in \mathcal{B}_x^0(z)$.  Let $\beta_\xi(u) = z$, so $J_{u}^{(r)} \beta_\xi = \eta \in \widehat{L}_\xi(z)$.  Now compare $J^{(r)} (\beta_\xi \circ \tau_u)$ and $J^{(r)} \beta_\eta$: both map the origin to $\eta$ and have image in $\widehat{L}_\xi \subset \mathcal{P}_x^{(r)}$.  Then $\beta_\xi \circ \tau_u = \beta_\eta$.  Next, $\eta = \omega \cdot h$ for a unique $h \in \mathcal{H}^{(r)}$.  Thus $\gamma^{-1} \beta = \beta_\omega^{-1} \beta_\xi = h \circ \tau_{-u}$. Now (4) is proved.  The same argument also shows that $\mathcal{B}^0(y)$ can be defined for all $y \in L_x$, assuming $x \in M_{\mathcal{R}}$.  This completes the proof of (1).

Given $\beta, \gamma$ as in point (2), one can write 
$$\mathcal{B}_x^0(y) = \beta \cdot \mathcal{H}^{(r),0} \qquad \mathcal{B}_x^0(z) = \gamma \cdot \mathcal{H}^{(r),0}$$
and the conclusion follows immediately.

Point (3) is just as in the proof of theorem \ref{thm.main} (1).

% $\widehat{L}_\xi \rightarrow L_x$-- surjective local diffeo.

%% Next consider $\beta_\xi, \beta_\eta \in \mathcal{B}_x$ with overlapping images in $L_x$.  Suppose $\beta_\xi({\bf 0}) = y$ and $\beta_\eta({\bf 0}) = z$.  Choose $w$ in the intersection of the images.

If the spectrum is $1/2$ pinched, then $r=1$, and $\mathcal{H}^{(1)}$ comprises linear transformations, so the atlas $\mathcal{B}_x^0$ is a flat affine structure on $L_x$. 
\end{Pf}

\begin{theorem}
\label{thm.foliation.connxns}
The foliation $\cup_x L_x$ contains a $Z(\varphi^t)$-invariant filtered family of $C^{q-1}$ subfoliations 
$$ L^{i_1} \subset \cdots L^{i_{l}} = L \qquad l \leq r,$$
a subfamily of the filtered family in theorem \ref{thm.ruelle.foliations}, equipped with flat connections $\nabla_x^j$ on the normal bundles of $L^{i_{j-1}} \subset L^{i_j}$ for all $1 \leq j \leq l$.  The $\oplus_x \nabla_x^j$ are $Z(\varphi^t)$-invariant, measurable in $x$, and $C^{q-2}$ inside almost every $L_x$; in particular, almost every leaf $L^{i_1}_x$ carries a $C^{q-2}$ invariant flat connection.
\end{theorem}

\begin{Pf}
Extract the components of the Lyapunov filtration corresponding to the conglomerated distributions:
$$V^{i_j} = \sum_{m \leq j} \widetilde{W}^m = \sum_{m \leq i_j} W^m$$ 
where $i_1, \cdots, i_{r_*} = i_l$ and the subspaces $\widetilde{W}^m$ are as in the proof of theorem \ref{thm.main} (2).  The $Z(F^t_x)$-invariant distributions $\mathcal{V}_x^j$, $1 \leq j \leq l$, on $\mathcal{E}_x$ are tangent at ${\bf 0}$ to $V^{i_j}$, and thus to the submanifolds $\nu_x^{i_j}$ given by theorem \ref{thm.ruelle.submanifolds}.  Then for almost all $x \in M$, the images $(\theta_x)_* \mathcal{V}_x^j$ are tangent to $L_x^{i_j}$, where $L^{i_j}$ is the foliation given by theorem \ref{thm.ruelle.foliations}.

The $Z(F^t_x)$-invariant flat connections $\nabla_x^j$ on $\mathcal{V}_x^j/\mathcal{V}_x^{j-1}$ given by theorem \ref{thm.main} (2) push forward under $\theta_x$ to flat connections on $T_xL^{i_j}/T_xL^{i_{j-1}}$, which are $Z(\varphi^t)$-invariant and $C^{q-2}$ inside $\mathcal{E}_x$.
\end{Pf}

\bibliographystyle{amsplain}
\bibliography{karinsrefs}

\begin{tabular}{l}
Karin Melnick   \\
Department of Mathematics  \\
4176 Campus Drive \\
University of Maryland  \\
College Park, MD 20742 \\
USA \\
karin@math.umd.edu 
\end{tabular}

\end{document}